\let\noarrow = t
\input eplain

\let\noarrow = t

\input eplain


\magnification=\magstep1

\topskip1truecm
\def\pagewidth#1{
  \hsize=#1
}

\def\pageheight#1{
  \vsize=#1
}

\pageheight{23.5truecm} \pagewidth{16truecm}

\abovedisplayskip=3mm \belowdisplayskip=3mm
\abovedisplayshortskip=0mm \belowdisplayshortskip=2mm
\parindent1pc

\normalbaselineskip=13pt \baselineskip=13pt

\voffset=0pc \hoffset=0pc


\newdimen\abstractmargin
\abstractmargin=3pc


\newdimen\footnotemargin
\footnotemargin=1pc


\font\eightrm=cmr8 \relax 
\font\sixrm=cmr6 \relax 
\font\eighti=cmmi8 \relax     \skewchar\eighti='177 
\font\sixi=cmmi6 \relax       \skewchar\sixi='177   
\font\eightsy=cmsy8 \relax    \skewchar\eightsy='60 
\font\sixsy=cmsy6 \relax      \skewchar\sixsy='60   
\font\eightbf=cmbx8 \relax 
\font\sixbf=cmbx6 \relax   
\font\eightit=cmti8 \relax 
\font\eightsl=cmsl8 \relax 
\font\eighttt=cmtt8 \relax 

\catcode`\@=11
\newskip\ttglue

\def\eightpoint{\def\rm{\fam0\eightrm}%
 \textfont0=\eightrm \scriptfont0=\sixrm
 \scriptscriptfont0=\fiverm
 \textfont1=\eighti \scriptfont1=\sixi
 \scriptscriptfont0=\fivei
 \textfont2=\eightsy \scriptfont2=\sixsy
 \scriptscriptfont2=\fivesy
 \textfont3=\tenex \scriptfont3=\tenex
 \scriptscriptfont3=\tenex
 \textfont\itfam\eightit \def\it{\fam\itfam\eightit}%
 \textfont\slfam\eightsl \def\sl{\fam\slfam\eightsl}%
 \textfont\ttfam\eighttt \def\tt{\fam\ttfam\eighttt}%
 \textfont\bffam\eightbf \scriptfont\bffam\sixbf
   \scriptscriptfont\bffam\fivebf \def\bf{\fam\bffam\eightit}%
 \tt \ttglue=.5em plus.25em minus.15em
 \normalbaselineskip=9pt
 \setbox\strutbox\hbox{\vrule height7pt depth3pt width0pt}%
 \let\sc=\sixrm \let\big=\eifgtbig \normalbaselines\rm}


 \font\titlefont=cmbx12 scaled\magstep1
 \font\sectionfont=cmbx12
 \font\ssectionfont=cmsl10
 \font\claimfont=cmsl10

 \font\normalfont=cmr10

\catcode`\@=11 \font\teneusm=eusm10 
\font\seveneusm=eusm7  \font\fiveeusm=eusm5
\newfam\eusmfam \textfont\eusmfam=\teneusm
\scriptfont\eusmfam=\seveneusm \scriptscriptfont\eusmfam=\fiveeusm
\def\hexnumber@#1{\ifcase#1
0\or1\or2\or3\or4\or5\or6\or7\or8\or9\or         A\or B\or C\or
D\or E\or F\fi } \edef\eusm@{\hexnumber@\eusmfam}
\def\euscr{\ifmmode\let\next\euscr@\else
\def\next{\errmessage{Use \string\euscr\space only in math mode}}\fi\next}
\def\euscr@#1{{\euscr@@{#1}}} \def\euscr@@#1{\fam\eusmfam#1} \catcode`\@=12

\catcode`\@=11 \font\teneuex=euex10 
 \font\seveneuex=euex7  \newfam\euexfam
\textfont\euexfam=\teneuex  \scriptfont\euexfam=\seveneuex
 \def\hexnumber@#1{\ifcase#1
0\or1\or2\or3\or4\or5\or6\or7\or8\or9\or         A\or B\or C\or
D\or E\or F\fi } \edef\euex@{\hexnumber@\euexfam}
\def\euscrex{\ifmmode\let\next\euscrex@\else
\def\next{\errmessage{Use \string\euscrex\space only in math mode}}\fi\next}
\def\euscrex@#1{{\euscrex@@{#1}}} \def\euscrex@@#1{\fam\euexfam#1}
\catcode`\@=12

\catcode`\@=11 \font\teneufb=eufb10 
\font\seveneufb=eufb7  \font\fiveeufb=eufb5
\newfam\eufbfam \textfont\eufbfam=\teneufb
\scriptfont\eufbfam=\seveneufb \scriptscriptfont\eufbfam=\fiveeufb
\def\hexnumber@#1{\ifcase#1
0\or1\or2\or3\or4\or5\or6\or7\or8\or9\or         A\or B\or C\or
D\or E\or F\fi } \edef\eufb@{\hexnumber@\eufbfam}
\def\euscrfb{\ifmmode\let\next\euscrfb@\else
\def\next{\errmessage{Use \string\euscrfb\space only in math mode}}\fi\next}
\def\euscrfb@#1{{\euscrfb@@{#1}}} \def\euscrfb@@#1{\fam\eufbfam#1}
\catcode`\@=12

\catcode`\@=11 \font\teneufm=eufm10 
\font\seveneufm=eufm7  \font\fiveeufm=eufm5
\newfam\eufmfam \textfont\eufmfam=\teneufm
\scriptfont\eufmfam=\seveneufm \scriptscriptfont\eufmfam=\fiveeufm
\def\hexnumber@#1{\ifcase#1
0\or1\or2\or3\or4\or5\or6\or7\or8\or9\or         A\or B\or C\or
D\or E\or F\fi } \edef\eufm@{\hexnumber@\eufmfam}
\def\euscrfm{\ifmmode\let\next\euscrfm@\else
\def\next{\errmessage{Use \string\euscrfm\space only in math mode}}\fi\next}
\def\euscrfm@#1{{\euscrfm@@{#1}}} \def\euscrfm@@#1{\fam\eufmfam#1}
\catcode`\@=12

\catcode`\@=11 \font\teneusb=eusb10 
\font\seveneusb=eusb7  \font\fiveeusb=eusb5
\newfam\eusbfam \textfont\eusbfam=\teneusb
\scriptfont\eusbfam=\seveneusb \scriptscriptfont\eusbfam=\fiveeusb
\def\hexnumber@#1{\ifcase#1
0\or1\or2\or3\or4\or5\or6\or7\or8\or9\or         A\or B\or C\or
D\or E\or F\fi } \edef\eusb@{\hexnumber@\eusbfam}
\def\euscrsb{\ifmmode\let\next\euscrsb@\else
\def\next{\errmessage{Use \string\euscrsb\space only in math mode}}\fi\next}
\def\euscrsb@#1{{\euscrsb@@{#1}}} \def\euscrsb@@#1{\fam\eusbfam#1}
\catcode`\@=12

\catcode`\@=11 \font\tenmsa=msam10 
\font\sevenmsa=msam7  \font\fivemsa=msam5
\font\tenmsb=msbm10  \font\sevenmsb=msbm7
 \font\fivemsb=msbm5 \newfam\msafam
\newfam\msbfam \textfont\msafam=\tenmsa
\scriptfont\msafam=\sevenmsa
  \scriptscriptfont\msafam=\fivemsa
\textfont\msbfam=\tenmsb  \scriptfont\msbfam=\sevenmsb
  \scriptscriptfont\msbfam=\fivemsb
\def\hexnumber@#1{\ifcase#1 0\or1\or2\or3\or4\or5\or6\or7\or8\or9\or
        A\or B\or C\or D\or E\or F\fi }
\edef\msa@{\hexnumber@\msafam} \edef\msb@{\hexnumber@\msbfam}
\mathchardef\square="0\msa@03 \mathchardef\subsetneq="3\msb@28
\mathchardef\supsetneq="3\msb@29 \mathchardef\ltimes="2\msb@6E
\mathchardef\rtimes="2\msb@6F \mathchardef\dabar="0\msa@39
\mathchardef\daright="0\msa@4B \mathchardef\daleft="0\msa@4C

\def\Bbb{\ifmmode\let\next\Bbb@\else
        \def\next{\errmessage{Use \string\Bbb\space only in math mode}}\fi\next}
\def\Bbb@#1{{\Bbb@@{#1}}}
\def\Bbb@@#1{\fam\msbfam#1}
\catcode`\@=12



\newcount\senu
\def\senum{\number\senu}
\newcount\ssnu
\def\ssnum{\number\ssnu}
\newcount\fonu
\def\fonum{\number\fonu}

\def\num{{\senum.\ssnum}}
\def\numfo{{\senum.\ssnum.\fonum}}


\outer\def\section#1\par{\vskip0pt
  plus.3\vsize\penalty20\vskip0pt
  plus-.3\vsize\bigskip\vskip\parskip
  \message{#1}\centerline{\sectionfont\senum\enspace#1.}
  \nobreak\smallskip}

\def\endsection{\advance\senu by1\penalty-20\smallskip\ssnu=1}
\outer\def\ssection#1\par{\bigskip
  \message{#1}{\noindent\bf\num\ssectionfont\enspace#1.\thinspace}
  \nobreak\normalfont}

\def\endssection{\advance\ssnu by1\smallskip\ifdim\lastskip<\medskipamount
\removelastskip\penalty55\medskip\fi\fonu=1\normalfont}

\def\proclaim #1\par{\bigskip
  \message{#1}{\noindent\bf\num\enspace#1.\thinspace}
  \nobreak\claimfont}

\def\cor{\proclaim Corollary\par}
\def\defi{\proclaim Definition\par}
\def\lemma{\proclaim Lemma\par}
\def\prop{\proclaim Proposition\par}
\def\rmk{\proclaim Remark\par\normalfont}
\def\thm{\proclaim Theorem\par}

\def\example{\proclaim Example\par\normalfont}

\def\endcor{\endssection}
\def\enddefi{\endssection}
\def\endlemma{\endssection}
\def\endprop{\endssection}
\def\endrmk{\endssection}
\def\endthm{\endssection}

\def\endexample{\endssection}

\def\Proof{{\noindent\sl Proof: \/}}


\def\maplefto#1{\ \smash{\mathop{\longleftarrow}\limits^{#1}}\ }

\def\llongrightarrow{\relbar\joinrel\relbar\joinrel\rightarrow}
\def\lllongrightarrow{\hbox to 40pt{\rightarrowfill}}

\def\twoheadrightarrow{\rightarrow\kern -8pt\rightarrow}

\def\maprighto#1{\smash{\mathop{\longrightarrow}\limits^{#1}}}

\def\llongmaprighto#1{\ \smash{\mathop{\llongrightarrow}\limits^{#1}}\ }

\def\lllongmaprighto#1{\ \smash{\mathop{\lllongrightarrow}\limits^{#1}}\ }

\def\longleftmapsto{\longleftarrow\kern-2pt\mapstochar\;}

\def\llongmapsto{\,\vert\kern-3.2pt\joinrel\longrightarrow\,}
\def\llongmapsto{\,\vert\kern-3.7pt\joinrel\llongrightarrow\,}
\def\lllongmapsto{\,\vert\kern-5.5pt\joinrel\lllongrightarrow\,}

\def\isomarrow{\maprighto{\lower3pt\hbox{$\scriptstyle\sim$}}}
\def\llongisomarrow{\llongmaprighto{\lower3pt\hbox{$\scriptstyle\sim$}}}
\def\lllongisomarrow{\lllongmaprighto{\lower3pt\hbox{$\scriptstyle\sim$}}}

\def\lisomarrow{\maplefto{\lower3pt\hbox{$\scriptstyle\sim$}}}

\font\labprffont=cmtt8
\def\strutdepth{\dp\strutbox}
\def\labtekst#1{\vtop to \strutdepth{\baselineskip\strutdepth\vss\llap{{\labprffont #1}}\null}}
\def\marglabel#1{\strut\vadjust{\kern-\strutdepth\labtekst{#1\ }}}

\def\label #1. #2\par{{\definexref{#1}{\num}{#2}}}
\def\labelf #1\par{{\definexref{#1}{\numfo}{formula}}}
\def\labelse #1\par{{\definexref{#1}{\num}{section}}}


\def\Spec{{\rm Spec}}

\def\Z{{\bf Z}}

\def\X{{\cal X}}
\def\O{{\cal O}}

\def\UX{\overline{X}}
\def\CC{{\Bbb C}}

\def\QQ{{\Bbb Q}}
\def\PP{{\Bbb P}}
\def\II{{\Bbb I}}
\def\MM{{\Bbb M}}
\def\AA{{\Bbb A}}
\def\ZZ{{\Bbb Z}}
\def\hX{\hat{X}}

\senu=1 \ssnu=1 \fonu=1

\centerline{\titlefont Horizontal sections of connections on
curves and transcendence}

\

\centerline{ C. Gasbarri}

\bigskip

\centerline{\today}

\bigskip
\bigskip

{\insert\footins{\leftskip\footnotemargin\rightskip\footnotemargin\noindent\eightpoint
$2000$ {\it Mathematics Subject Classification}. Primary 11J91,
14G40, 30D35
}

\vbox{{\leftskip\abstractmargin \rightskip\abstractmargin
\eightpoint

\noindent A{{\sixrm BSTRACT}}: Let $K$ be a number field, $\UX$ be
a smooth projective curve over it and $D$ be a reduced divisor on
$\UX$. Let $(E,\nabla)$ be a fibre bundle with connection having
meromorphic poles on $D$. Let $p_1,\dots,p_s\in\UX(K)$ and
$X:=\UX\setminus\{D,p_1,\dots, p_s\}$ (the $p_j$'s may be in the support of
$D$). Using
tools from Nevanlinna theory and formal geometry, we give the
definition of $E$--section of type $\alpha$ of the vector bundle
$E$ with respect to the points $p_j$; this is the natural generalization of the
notion of $E$
function defined in Siegel Shidlowski theory. We prove that the
value of a $E$--section of type $\alpha$ in an algebraic point different from
the $p_j$'s has
maximal transcendence degree. Siegel Shidlowski theorem is a
special case of the theorem proved. We give an application to
isomonodromic connections.

}}

\

\section Introduction\par

\

Many questions in transcendental theory may be resumed in this
"meta--question": Suppose that $U$ is a variety defined over a
number field $K$ and $G(F, F^{(1)},\dots, F^{(n)})=0$ is an
algebraic system of differential equations defined over $U$ (the
functions defining $G$ are in $K(U)$). Suppose that $F:=(F_1,\dots
F_n)$ is a local solution of the differential equation. Let $q\in
U(K)$; what can we say about $Trdeg_{\QQ}(K(F(q)))$? Up to the
fact that this degree of transcendency is  bounded from above  by
$Trdeg_{K(U)}(K(U)(F))$, we cannot say a lot about this question
in general.

If we restrict our attention to {\it systems of linear
differential equations} over the projective line and regular over
the multiplicative group $\Bbb G_m$, the Siegel Shidlowski theory
give us a very powerful and satisfactory answer. Let's recall the
main result of the theory (in a simplified version), cf. for
instance [La]:

Let \labelf siegelslisys\par$${{dY}\over{dz}}=AY\;\;\;\;\;\;{\rm
with}\;\;\;\; A\in M_n(\QQ(z))\eqno{{(\numfo)}}$$ be a linear
system of differential equations. Suppose that $F=(f_1(z),\dots,
f_n(z))$ is a solution of \ref{siegelslisys} with the following
properties:

(a) the functions $f_1(z),\dots, f_n(z)$ are algebraically
independent over $\CC(z)$;

(b) Each of the $f_i(z)$ has a Taylor expansion
$f_i(z)=\sum_{j=0}^{\infty}a_{ij}{{z^j}\over{j!}}$ with
$$a_{ij}\in\QQ \;\;\;\;\;\;{\rm and}\;\;\;\;\;\; H(a_{ij})\ll
j^{\epsilon j}$$ ($H(\cdot)$ being the exponential height).

Then, for every $q\in\QQ^\ast$, we have that
$Trdeg_\QQ(\QQ(f_1(q),\dots, f_n(q))=n$. Recall that functions
with  property (b) above are called $E$--functions.

It is well known that the criterion above has many important
consequences; in particular, the Hermite Lindeman  Theorem is a
special case of it (take $f_i(z)=e^{\alpha_iz}$) and  non trivial
transcendence properties of special values of some hypergeometric
and Bessel functions can be deduced from it.

Of course, if one could generalize Siegel Shidlovski theory to
arbitrary variety, the general "meta--question" above would have a
satisfactory answer in the linear case. Unfortunately, the example
below shows that a believable statement over an arbitrary variety
is not easy to find:

\example Let $f_1(z_1,z_2)=e^{z_1}$, $f_2(z_1,z_2)=e^{z_2}$ and
$f_3(z_1,z_2)=z_1$; then $f_i$ are algebraically independent over
$\CC(z_1,z_2)$, they satisfy a system of linear differential
equations with coefficients in $\QQ[z_1,z_2]$ but, for every $a\in\QQ$ the
restriction of them to the line $z_1=az_2$ are algebraically
dependent. Perhaps one have to look for interactions between
Siegel Shidlowski theory with the theorem of Ax [Ax].
\endexample

In this paper we develop a transcendental theory, analogous to the
Siegel--Shidlovski's, for horizontal sections of connections over
arbitrary curves. As the example below shows, even in this case
some caveat are necessary:

\label ordertwo. example\par\example consider the functions
$f_1(z)=e^z$ and $f_2(z)=e^{z^2}$. They are algebraically
independent, they satisfy a system of differential equation with
coefficients in $\QQ[z]$, but for every $a\in\QQ$, we have that
$f_1(a)$ and $f_2(a)$ are algebraically dependent.
\endexample

This tells us that a condition analogous to (b) is necessary.
Observe  that an $E$--function is an entire function with order of
growth one! In general, the lemma below shows that horizontal sections of
vector bundles with connections are of finite growth, for a precise definition
of finite order of growth, cf. \S 4:

\label orderofgrown1. proposition\par\prop Let $\Delta$ be the
unit disk with coordinate $z$.  Let
$$Y'={{A}\over {z^n}}\cdot Y$$
with $A\in\O_\Delta$ be a system of linear differential equations
with poles only in the origin of order at most $n$. Let $F$ be an
analytic solution in $\Delta^\ast:=\Delta\setminus\{ 0\}$. Then
$$\log\Vert F\Vert\leq {{C}\over{\vert z\vert^{n-1}}}$$ for a
suitable positive constant $C$.

\endprop

\Proof The statement is standard, we give a sketch of proof for reader
convenience: Let $\vert z\vert=1$ and $t\in (0,1)$; consider the
function $h(t):=\Vert F(tz)\Vert$. Then ${{d}\over{dt}}h^2(t)\leq
2\vert F(tz)F'(tz)\vert$, consequently $\vert h'(tz)\vert \leq
\Vert F'(tz)\Vert$. This implies that there is a constant $C$ such
that
$$\left\vert {{h'(t)}\over {h(t)}}\right\vert\leq {{C}\over{t^n}}.$$
The conclusion follows integrating both sides with respect to $t$.

\smallskip

On the other direction, in our paper [Ga1], as a consequence of the main results
proved there, we can deduce the following:

\label schneiderlang. theorem\par\thm Let $\overline X/_\QQ$ be a
smooth projective curve and $D\subseteq\UX$ be a reduced divisor.
Denote by $X$ the affine curve $\UX\setminus D$ Let $(E;\nabla)$
be a fibre bundle with connections having meromorphic poles along
$D$. Let $f_\CC:X(\CC)\to E$ be an horizontal section of finite
order of growth  $\rho$. Then
$$Card(f(X(\CC))\cap E(\QQ))\leq
{{rk(E)+2}\over{rk(E)}}\cdot\rho.$$
\endthm
(The theorem above is not explicitly stated in [Ga1], it can
obtained as a particular case of theorem 1.1 of loc. cit.). A
linear algebra manipulation may improve \ref{schneiderlang} to
obtain:

\label corollaofsl. corollary\par\cor In the hypotheses of theorem
\ref{schneiderlang} we have that
$$Card(f(X(\CC))\cap E(\overline{\QQ}))\leq \rho.$$
\endcor

{\it Sketch of Proof}: Of course, if we apply \ref{schneiderlang}
to the symmetric power of $E$ with the induced connection and the
induced section we obtain that $Card(f(X(\CC))\cap E(\QQ))\leq
\rho$.

Since $f$ is the horizontal section of a section, it defines a
$LG$ germ of type 1 (definition in [Ga \S 3]) on every rational
point of $X$. Let $p_1,\dots , p_r$ be points such that $f(p_j)\in
E(\overline\QQ)$.  If we multiply a $LG$-- germ of type with an
integral section (over $Spec(\ZZ)$) of a line bundle we obtain
again an $LG$ germ and also the order of growth do not change.
Thus we may suppose that the natural restriction map
$H^0(E^\vee)\to\oplus_jE^\vee\vert_{p_j}$ is surjective. Fixing a
trivialization of $E$ at each $p_j$, we see that we can find two
sufficiently generic  polynomials $H^1_j$ and $H^2_j$ such that
$H^i_j(f\vert_{p_j})=0$. Because of the condition above, we can
construct two global sections $H_i\in
\oplus_{h=0}^nSym^h(E^\vee)$, for a suitable $n$,  such that
$H_i\vert_{p_j}=H^i_j$. Consequently $F_H:=(H_1(f),H_2(f))$ is an
analytic section of $\O_X\oplus\O_X$ which has  order of growth
$\rho$, vanish at the $p_j$ and it is an  $LG$--germs in these
points. Due to the fact that $f$ is Zariski dense and that the
$H^i_j$ are generic, one verify that $F_H$ is Zariski dense. Thus
one can apply the first part of the proof and conclude.

\smallskip

Thus the corollary above claims that, under the condition
that the order of growth is $\rho$, there are at most $\rho$ algebraic
point on the image of $f$; a small variation of the argument shows that, over
an arbitrary number field $K$, there can be at most $\rho[K:\QQ]$
algebraic values. Moreover, observe that, the functions
of example \ref{ordertwo} are of total order two and it is not difficut to
deduce from Schneider--Lang theorem that there is only one point ($z=0$) where
the value of both functions is algebraic. Thus we see that, if an horizontal
section of a vector bundle has less algebraic values on algebraic points then
the order of growth, we cannot say anything about the algebraic independence of
the values of  the section on other algebraic points. In the limit case
when the number of algebraic values is the same then the order of growth, we
can say more.

This brings to the definition of $E$--section of type $\alpha$ of a
vector bundle. Roughly speaking, $X$ be an affine curve defined
over $\QQ$, $V$ be vector bundle over it, $\alpha$ is an integer
(this is not necessary in general but here we suppose for
simplicity) and and $p_1,\dots, p_s\in X(\QQ)$; an analytic section
$f:X(\CC)\to E$ is said to be an $E$--section of type $\alpha$, if
locally around each of the $p_j$ we may write $f=(f_1(z),\dots, f_m(z))$ with
$$f_i(z)=\sum_{j=0}^\infty
a_{ij}{{z^j}\over{(j!)^\alpha}}\;\;\;\;\;\; {\rm with}\;\;\;\;\;
a_{ij}\in \ZZ[{{1}\over{N}}]$$ and the order of growth of $f$ is
${{s}\over{\alpha}}$. The order of growth is defined using
Nevanlinna theory on the Riemann surface $X(\CC)$. The precise
definition is more involved and is given over an arbitrary number
field and arbitrary $\alpha$, cf \S 6. The $E$--sections of type
$\alpha$ are a good generalization of $E$--functions over
arbitrary curves. Never the less it is important to notice that
while the local behavior and the growth behavior of an $E$
function is resumed in its definition as a power series, the local
and global properties of an $E$ sections are defined separately
via formal geometry and Nevanlinna theory. In the paper [Be1] the
author proves a generalization of the Schneider--Lang criterion
just imposing a local Gevrey condition which is very similar to
our definition of $LG$--germs. With this definition in mind we can
state the main theorem proved in this paper (here, for simplicity,
we state it over $\QQ$, for the general statement cf.
\ref{SSoncurves}):

\label SSoncurvesintro. theorem\par\thm Let $\UX/_\QQ$ be a smooth
projective curve. Let $D$ be an effective divisor on $\UX$ and
$(E,\nabla)$ be a fibre bundle of rank $m>1$ with connection with
meromorphic poles along $D$. Let $p_1,\dots,p_s\in \UX(\QQ)$ be  rational
points. Let $D':=D-\{p_1,\dots, p_s\}$ and $X:=\UX\setminus D'$. Let
$f:X(\CC)\to E$
be an analytic  {\rm horizontal section} which is an $E$--section
of type $\alpha$ with respect to the points  $p_j$. Suppose that the image of
$f$ is Zariski dense in $E$. Let $q\in X(\QQ)$, then
$$Trdeg_\QQ(\QQ(f(q)))=m.$$
\endthm

Observe that if $X=\PP^1$, $D=0+\infty$, and we have only one point $p=0$ we
find the classical theorem by Siegel and Shidlowski. The requirement that the
image is Zariski dense is equivalent to the requirement that the entries
of $f$ are algebraically independent over $\overline\QQ(\overline
X)$.

Even in the case when $\UX=\PP^1$ but $D$ is arbitrary, theorem
\ref{SSoncurvesintro} is stronger then the classical theorem by
Siegel and Shidlowski. This is due to the fact that the recent
papers [A1], [A2] and [Beu] prove the following: suppose we have a
system of differential equations as in \ref{siegelslisys} and a
solution $F:=(f_1,\dots ,f_n)$ such that $f_i$ are $E$--functions.
Then there is a system of differential equations $Z'=GZ$ with
$G\in M_n(\QQ[z;{{1}\over{z}}])$, a solution $E$ of it and a
matrix $M\in M_n(\overline\QQ[z])$ such that $F=M\cdot E$. Thus if
we restrict our attention to $E$--functions, or to functions which
are defined by power series having similar conditions, all the
transcendency results can be deduced from the classical
Siegel--Shidlowki statement.

Roughly speaking the statement of \ref{SSoncurvesintro} may be interpreted in
the following way: Suppose that we have a solution $f$ of a linear differential
equation over an affine curve and the order of growth is $\rho$. If $f$ has
algebraic value on exacly $\rho$ rational  points, then the values of the
entries of $f$ on any other rational point are algebraically independent. 
Example \ref{ordertwo} tells us that we cannot hope better then this.

Some words on the methodology. It is well known that there is an
analogy between the arithmetic of varieties over number fields,
the arithmetic of varieties over function fields and Nevanlinna
theory (cf. for instance [Vo]). This analogy is a dictionary which
allows to translate statements from a theory to another. In this
paper (and in our previous [Ga1]) this philosophy is pushed
forward: Instead of just an analogy, we try to develop an unified
theory where arithmetic, analysis and algebraic geometry interact
together. Usually we can find non trivial upper bounds from
Nevanlinna theory and lower bounds from arithmetics, Algebraic
geometry gives a common language where these two tools may
interact. The papers [Be2] and [Be3] are, in our opinion, very
near to the spirit of this paper; on these papers also, deep
generalizations of Siegel Shidlowki theory are given.

This paper is organized as follows. In \S 2 we give some general
criteria for an element of a complex vector space  to have
algebraically independent coordinates. In \S 3 we prove a zero
lemma over arbitrary curve which replace the classical Shidlowski
Lemma; the statement is formally the same then the Shidlowski
Lemma, but the proof is simpler and use some tools from algebraic
geometry: theory of vector bundles and Hilbert schemes. In \S 4 we
explain the tools from Nevanlinna theory which are needed; we use
a version of Nevanlinna theory (developed in [Ga1]) which allows
to prove powerful lemmas  of Scwhartz's type over (special kind
of) Riemann surfaces. In \S 5 and 6 we develop the notion of $E$
sections of type $\alpha$ and explain the main properties of them.
Eventually in \S 7 we state and prove the main theorem of the
paper.

\ssection An application\par We can give an interesting
application  to the theory of the isomonodromic connections. Let
$X$ be a curve and $(E_1;\nabla_1)$ and $(E_2,\nabla_2)$ two
integrable connections of rank $n$ over it. Let
$\rho_i:\pi_1(X)\to GL_n$ be the monodromy representation
associated to $(E_i,\nabla_i)$. Suppose that $\rho_1$ is
equivalent to $\rho_2$; thus the trivial representation is a
subrepresentation of $\rho_1\otimes\rho_2^\vee$; consequently we
get a global horizontal section of $E_1\otimes E_2^\vee$. Provided
that it has the right order of grown, this section is the typical
section to which we can apply the criterion.

We may guarantee the right order of grown by proposition \ref{orderofgrown1}. 
In particular it guarantees  that if we have a connection on a
projective curve with poles at most of order two, then an
horizontal section will define an $E$--section of type one (cf.
definition   \ref{esection}) over any smooth point rational point
of the connection.

From this we can obtain the following: Let $X$ be a projective
curve over $\Bbb Q$. Let $D$ be a reduced effective divisor over
$X$. Denote by $U$ the affine curve $X\setminus D$. Let
$(E_1;\nabla_1)$ and $(E_2,\nabla_2)$ be two fibre bundle with
connection having poles at most of order two along $D$. Suppose
that the corresponding representations $\rho_i:\pi_1(U)\to GL_N$
are isomorphic. Let $q\in U(\Bbb Q)$. We can find an analytic
isomorphism $\varphi:E_1\to E_2$ over $U$, which restrict to the
identity over $q$. Suppose that $Trdeg_{\QQ(X)}(\QQ(\varphi))=r$

\label isomomodromic2. theorem\par\thm Let $V$ be an analytic
neighborhood of $q$ and $p\in V\cap U(\Bbb Q)$. Let $F$ be a local
horizontal section around $q$ of $(E,\nabla_1)$. Suppose that
$F(p)\in\overline{\Bbb Q}$, then
$Trdeg_{\QQ}(\QQ(\varphi(F(p))=r$.

\endthm

The proof is a direct application of Theorem \ref{SSoncurves}.

A non trivial way to construct examples where theorem
\ref{isomomodromic2} apply, is the following: Let $B$ be a reduced
divisor in $\AA^1_\QQ$. Let $X$ be a smooth projective curve
defined over $\QQ$. Let $D$ be a reduced divisor over $X$. Over
$X\times\AA^1$ consider the divisors $H_1= D\times \AA^1$ and
$H_2=X\times B$ and $H=H_1+H_2$. Suppose that $(E,\nabla)$ is a
fibre bundle with {\it integrable} connection over $X\times\AA^1$
with poles around $H$ and which are at most of order two around
$H_1$.

Then, for every $q\in\AA^1(\QQ)\setminus B$ the restriction
$(E_q,\nabla_q)$ of $(E,\nabla)$ to $X\times\{ q\}$ is a vector
bundle with integrable connection having poles of order at most
two along $D$.

By construction, for every couple $q_i$, $q_2\in\AA^1(\QQ)$, the
vector bundles $(E_{q_i},\nabla_{q_i})$ have conjugated monodromy.
Thus the theorem apply in this case.

\endssection

\ssection An explicit example\par  Let $a, b, c\in \QQ$, for every
$x\in\QQ$ consider the linear system of differential equations
$$\nabla_x :\;\;\;\;\;\; {{dY}\over{dz}}=\left(
{{1}\over{z^2}}\cdot\left(\matrix{a& (a-b)x\cr 0 &
b\cr}\right)+{{1}\over{z}}\cdot\left(\matrix{(1-x) & -x^2\cr 1 & 1
\cr}\right)  +{{1}\over{z-1}}\cdot\left(\matrix{1 & c\cr 0 &
1\cr}\right)\right)\cdot Y.$$ Then,  up to conjugation, for every
couple $x_0, x_1\in \Bbb Q$ the linear systems $\nabla_{x_0}$ and
$\nabla_{x_1}$ have the same monodromy.

\Proof Consider, over $\PP^1\times\AA^1$, with local coordinates
$(z,x)$,  Denote by $\omega$ the matrix of differential forms
$$\omega:=\left(
{{1}\over{z^2}}\cdot A(x)+{{1}\over{z}}\cdot B(x)
+{{1}\over{z-1}}\cdot\left(\matrix{1 & c\cr 0 &
1\cr}\right)\right)dz+{{1}\over{x}}\cdot\left(\matrix{1 & 1\cr 0
&1}\right)dx.$$ With $A(x)$ and $B(x)$ unknown matrices to be
determined. The system of differential equations
$${\cal E}:\;\; \nabla(Y)=\omega\cdot Y$$ defines a fibre bundle
with integrable connection if and only if
$$d(\omega)=\omega\wedge\omega.$$ Thus ${\cal E}$ is integrable if
and only if $A(x)$ and $B(x)$ are solution of the linear
differential system \labelf differentialsys\par$${{d
W(x)}\over{dx}}={{\left[ W(x);\left(\matrix{ 1 & 1\cr 0 &
1\cr}\right)\right]}\over{x}}.\eqno{{(\numfo)}}$$

A basis of solutions of the system \ref{differentialsys} is
$$\left\{\left(\matrix{1 &\log(x)\cr 0& 0\cr}\right)\; ;\; \left(\matrix{0 & 1\cr 0& 0\cr}\right)\; ;\;
\left(\matrix{-\log(x) &-\log^2(x)\cr 1& 0\cr}\right)\;
;\;\left(\matrix{0 &-\log(x)\cr 0& 1\cr}\right)\right\}$$ thus if
we put $\tilde x=\log(x)$ and choose $\nabla_0$ to be
$$\nabla_0 :\;\;\;\;\;\; {{dY}\over{dz}}=\left(
{{1}\over{z^2}}\cdot\left(\matrix{a& 0\cr 0 &
b\cr}\right)+{{1}\over{z}}\cdot\left(\matrix{1 & 0\cr 1 & 1
\cr}\right)  +{{1}\over{z-1}}\cdot\left(\matrix{1 & c\cr 0 &
1\cr}\right)\right)\cdot Y$$ the conclusion follows.
\smallskip

\rmk One can see that the method we used above leave a lot of
freedom in the choices; one can use all the powerful theory of
isomonodromic deformations of connection developed for instance in
[Ma]. Many cases that one may construct from this method give rise
to example which, by a little of work, may deduced from the
classical theorem of Siegel and Shidlowski. We were not able to
reduce the transcendency problem arising from the family of
differential equations proposed above from the classical theory.
Unfortunately we are not able to prove that the analytic
isomorphism between two different elements of the family has an
high degree of transcendency over $\QQ(z)$. In any case, it
doesn't seem to us that we can obtain such an isomorphism from
exponentials and algebraic functions.
\endrmk

\endssection

\endssection

\endsection

\section Criteria for algebraic independence\par

\

Let $K$ be a number field and $O_K$ be its ring of integers. We
fix an embedding $\sigma : K\to \CC$. Let $E$ be an hermitian
$O_K$ module of rank $m$. If $V$ is an hermitian $O_K$ module,
denote by $V_K$ the $K$ vector space $V\otimes_{O_K}K$ and by
$V_\CC$ the $\CC$ vector space $V\otimes\CC$ ($\CC$ is an $O_K$
module via $\sigma$). In this section we want to describe some
criteria that imply that the coordinates of an element $f\in
E_\CC$ are algebraically independent over $K$.

In the sequel, for every integer $n$, we denote by $E_n$ the
hermitian $O_K$ module $Sym^n(E)$.

Let $f\in E_\CC$. Denote by $f_n$ the image of $f^{\otimes n}$ in
$(E_n)_\CC$. For every $n$ denote by $V_n$ the smallest
$K$--subspace containing $f_n$ and by $r_n$ its dimension.

\label algebraicallyfree. definition\par\defi We will say that $f$
is algebraically independent over $K$ if $V_n=E_n$ for every
positive integer $n$.
\enddefi

\rmk If we fix a basis of $E_K$, $f$ is algebraically independent
over $K$ if the coordinates of it are transcendental numbers
algebraically independent over $\QQ$.\endrmk

We will now give some criteria which imply that $f$ is
algebraically independent over $K$. First of all we observe the
following trivial fact:

-- Suppose that $V_1$ and $V_2$ are vector spaces and
$\dim(V_2)<\dim(V_1)$  then
$$\lim_{n\to
\infty}{{\dim(Sym^n(V_2))}\over{\dim(Sym^n(V_1))}}=0.$$

The proof is trivial and left to the reader.

\label linearindependencecriterion. lemma\par\lemma The vector $f$
is algebraically independent over $K$ if and only if there is a
constant $c>0$ such that for every $n$,
$${{\dim(V_n)}\over{\dim (E_n)}}\geq c.$$
\endlemma
\Proof If $f$ is algebraically independent over $K$ then, by
definition we may take $c=1$.

Conversely, suppose that $f$ is {\it not} algebraically
independent over $K$, then there is an $n$ and a non trivial
subspace $V_n\subsetneq E_n$ containing $f_n$.  Thus for every
integer $m$, $f_{nm}\in Sym^m(V_n)\subsetneq E_{nm}$.
Consequently, there is a subsequence $n_m$ such that
$$\lim_{m\to\infty}{{\dim(V_{n_m})}\over{\dim(E_{n_m})}}=0.$$
The conclusion follows.

\smallskip

Observe that, the constant $c$ of the lemma above, either is zero
or, a posteriori, it is one.

We will give now a criterion which implies the hypotheses of lemma
\ref{linearindependencecriterion}. Before it we need to recall the
definition of the Arakelov degree.

-- If $M$ is an hermitian line bundle over $\Spec(O_K)$ we will
define its Arakelov degree by the following formula: Let $s\in
M\setminus\{ 0\}$; then
$$\widehat{\deg}(M):=\log({\rm Card}(M/s\cdot O_K))-\sum_{\sigma\in
M_{\infty}}\log\Vert s\Vert_\sigma.$$ This formula is well defined
because of the product formula (cf. for instance [SZ]).

-- If $\overline E$ is an arbitrary hermitian vector bundle over
$\Spec (O_K)$ then the line bundle $\bigwedge^{max}E$ is
canonically equipped with an hermitian metric; consequently we can
define the hermitian line bundle $\bigwedge^{\max}\overline E$. We
then define $\widehat{\deg}(\overline
E):=\deg(\bigwedge^{\max}(\overline E))$.

-- (Cramer rule): If $E$ is an hermitian vector bundle of rank
$r$, then there is a canonical hermitian isomorphism of vector
bundles:
$$\det(E)\otimes E^\vee\simeq \bigwedge^{r-1}E.$$

\smallskip

In the following we will denote by $c_i$'s constants which are
independent on  the rank $m$. Since, due to lemma
\ref{linearindependencecriterion}, in order to prove algebraic
independence, we can work with symmetric products, we may always
suppose that each $c_i$'s is very small compared with $m$.

We fix an hermitian line bundle $H$ over $\Spec(O_K)$. In the
sequel, if $F$ is an hermitian $O_K$ module, for every integer $x$,  we will denote $F(x)$ the
hermitian vector bundle $F\otimes H^{\otimes x}$. For $P_i\in
E^\vee(x)$ we denote by $F_i$ the vector $<P_i;f>\in H^{\otimes
x}_\CC$.

\label linearindependencecriterion2. theorem\par\thm Suppose that
there exist constants $c_i$ and $b_j$ (the former, possibly,
depending on $m$) for which the following holds: For every $x$
sufficiently big, there exist $P_1,\dots , P_m\in E^\vee(x)$ such
that:

-- They are linearly independent.

-- $\sup_{\sigma\in M_{\infty}}\{\log\Vert P_i\Vert_\sigma\}\leq
c_1x\log(x)+b_1x$.

-- $\sup_{\sigma\in M_{\infty}}\{\log\Vert F_i\Vert_\sigma\}\leq
c_1x\log(x)-c_2\cdot m\cdot x\cdot\log(x)+c_3\cdot x\cdot
\log(x)+b_2x$.

Then there are constants $C_i$ depending only on the $c_i$'s (thus
independent on $m$) such that
$$r_1\geq C_1m+C_2.$$

\endthm

\Proof Denote by $V_K\hookrightarrow E_K$ the minimal
$K$--subspace containing $f$. Let $V:=V_K\cap E$, then
$r_1=rk(V)$. For every positive integer $x$, denote by
$\tilde{P_i}$ the image of $P_i$ in $V^\vee(x)$. Observe that
there are constants $d_j$ such that
$\widehat{\deg}(V^\vee(x))=d_1+d_2\cdot x$. We can find $r_1$
elements within the $\tilde{P_i}$ which are linearly independent.
without loss of generality, we may suppose that they are
$\tilde{P_1},\dots , \tilde{P_{r_1}}$. The isomorphism of Cramer
rule give rise to the following equality
$$(\tilde{P_1}\wedge\dots\wedge\tilde{P_{r_1}})\otimes
f=\sum_i(-1)^i(\tilde{P_1}\wedge\dots\wedge\hat{\tilde{P_i}}\wedge\dots\wedge\tilde{P_{r_1}})\otimes
F_i.$$ Since $\tilde{P_1}\wedge\dots\wedge\tilde{P_{r_1}}$ is an
{\it integral} section of $V^\vee(x)$, then
$$\log\Vert\tilde{P_1}\wedge\dots\wedge\tilde{P_{r_1}}\Vert_\sigma\geq
d_1+d_2\cdot x-([K:\QQ]-1)(r_1\cdot c_1\cdot x\log(x)+b_1\cdot
x).$$ Thus we find
$$\eqalign{&d_1+d_2\cdot x-([K:\QQ]-1)(r_1\cdot c_1\cdot x\log(x)+b_1\cdot
x)+d_3\leq\cr &\leq (r_1-1)c_1\cdot x\log(x)+ c_1\cdot
x\log(x)-m\cdot c_2\cdot x\log(x)+c_3\cdot x\log(x)+b_2\cdot x+
d_4;\cr}$$ where the constants $b_i$'s, $c_i$'s and $d_i$'s are
independent on $x$. We divide everything by $x\log(x)$ and let $x$
tend to infinity and obtain
$$r_1\cdot c_1[K:\QQ]-m\cdot c_2+c_3\geq 0.$$ The conclusion
follows.

As corollary of theorem \ref{linearindependencecriterion2} and
lemma \ref{linearindependencecriterion} we find:

\label linearindependencecriterion3. corollary\par\cor Suppose
that, for every $n$ sufficiently big, we may apply theorem
\ref{linearindependencecriterion2} to $f_n$ in $(E_n)_\CC$,  with
the  $c_i$'s independent on $n$, then $f$ is algebraically
independent over $K$.
\endcor

\endsection

\

\section Connections and the Zero Lemma\par

\

Let $K$ be a number field embedded in $\CC$ and $\overline K$ be
its algebraic closure in $\Bbb C$. Let $X_K$ be a smooth
projective curve over $K$ and $D:=\sum_in_iP_i$ be an effective
divisor on $X_K$. We will denote by $F$ the function field in one
variable $\CC(X_K)$.

If $\nabla_i:F_i\to F_i\otimes\Omega^1_{X_K}(D_i)$ ($i=1,2$) are
fibre bundles with singular connections on $X_K$, then the tensor
product $F_1\otimes F_2$ is naturally equipped with a singular
connection $\nabla_{1,2}:F_1\otimes F_2\to F_1\otimes
F_2\otimes\Omega^1_{X_K}(l.c.m.(D_1,D_2))$ (where if
$D_i:=\sum_jn_{i,j}P_j$, then
$l.c.m.(D_1,D_2):=\sum_i\max_j\{n_{i,j}\}P_j)$.

The standard derivation $d:\O_{X_K}\to \Omega^1_{X_K}$ induces,
for every point $P\in X_K(\overline K)$, a singular connection
$\nabla^P:\O_{X_K}(P)\to \O_{X_K}(P)\otimes\Omega^1_{X_K}(P)$.
Thus, for every divisor $D$ the line bundle $\O(D)$ is equipped
with a canonical connection
$\nabla^D:\O(D)\to\O(D)\otimes\Omega^1_{X_K}(\vert D\vert)$
(where, if $D:=\sum_in_iP_i$ is a divisor, we define $\vert
D\vert$ to be the divisor $\sum_i\min_i\{1; \vert n_i\vert\}P_j$).

In the following we will denote by $H$ the line bundle
$\O_{X_K}(D)$; observe that, by the construction above, $H$ is
canonically equipped with a singular connection $\nabla^H:H\to
H\otimes\Omega^1_{X_K}(D)$; let $s\in H^0({X_K},H)$ be a section
such that $div(s)=D$. If $F$ is a coherent sheaf on ${X_K}$ and
$x$ an integer, we will denote by $F(x)$ the sheaf $F\otimes
H^{\otimes x}$.

Let $(E;\nabla^E)$ be a vector bundle on ${X_K}$ of rank $m$ with
a singular connection
$$\nabla^E : E\longrightarrow E\otimes\Omega^1_{{X_K}}(D).$$
For every integer $x$, the vector bundle $E(x)$ is equipped with a
singular connection $\nabla^x:E(x)\to
E(x)\otimes\Omega^1_{X_K}(D)$.

Fix a point $Q\in {X_K}(\CC)$, which {\it may be in the support of
$D$}. Let $\partial$ be a global section of $T_{X_K}(D)$ which do
not vanish on $Q$. We fix a section $s'\in H^0({X_K};H)$ such that
$s'(Q)\neq 0$. If necessary, we may enlarge $D$, thus we may
suppose that such a section exists as soon as we suppose that $D$
is of sufficiently big degree. The derivation $\partial$ and the
connection $\nabla^x$ induce a derivation
$$\nabla_\partial^x: E(x)\longrightarrow E(x+2).$$

Denote by $\hat X_Q$ the completion of ${X_K}$ around $Q$ and
denote by $E_Q$, $\Bbb L_Q$ etc. the restriction of $E$, $\Bbb L$
etc. to $X_Q$. We denote by $F_Q$ the completion of $F$ with
respect to the discrete valuation induced by $Q$.

The restriction of $(E(x);\nabla^x)$ to the generic fibre is a
$\cal D$--module $(E_F;\nabla^{x})$.

\label genericconnection. lemma\par\lemma Let $G\hookrightarrow
E(x)$ be a sub bundle (the quotient is without torsion). The
following properties are equivalent:

a) The bundle $G$ is a sub bundle with connection: the image of
$G$ via $\nabla^x$ is contained in $G\otimes\Omega^1_{X_K}(D)$.

b) The $F$ vector space $G_F$ is a sub $\cal D$--module of
$E(x)_F$: $\nabla(G_F)$ is contained in
$(G\otimes\Omega^1_{X_K}(D))_F\subseteq(E(x)\otimes\Omega^1_{X_K}(D))_F$.

c) The image of the $F$ vector space $G_F$ under the map
$\nabla^x_\partial$ is contained in $G(2)_F\subseteq E(x+2)_F$.
\endlemma

The proof is left to the reader.

Let $P\in H^0({X_K},E(x))$;  denote $P=P_0$  and
$P_{i+1}:=\nabla^{x+2i}_\partial(P_i)$; fix a positive integer
$r\leq m$ and denote by $\tilde P_i$ the elements $P_i\otimes
(s')^{\otimes 2(r-i)}$. The elements $\tilde P_0,\dots,\tilde P_r$
are elements of $H^0({X_K};E(x+2r))$. Let $G\subseteq E(x+2r)$ be
the sub vector bundle generated by the $\tilde P_i$ (the quotient
is without torsion). From the lemma above we deduce

\label connectionclosure. lemma\par\lemma Suppose that the $\tilde
P_i$ are linearly dependent as elements of $E(x+2r)_F$ then $G$ is
a sub bundle with connection of $E(x+2r)$.
\endlemma

Given a global section $P\in H^0({X_K},E(x))$, Suppose that
$\tilde P_0,\dots,\tilde P_r$ are linearly independent over $F$
and $\tilde P_0,\dots,\tilde P_{r+1}$ linearly dependent; denote
by $G$ the sub vector bundle generated by the $\tilde P_i$'s.

\label minimalsubbundle. definition\par\defi Given a global
section $P\in E(x)$, we will call the sub bundle $G\hookrightarrow
E(x+2r)$ constructed above, {\rm the minimal sub bundle with
connection generated by $P$}.
\enddefi

In particular we remark that if $\tilde P_0,\dots,\tilde P_{m-1}$
are linearly independent over $F$, then $G=E(x+2m)$.

Let $f$ be an horizontal section of $E^\vee_Q$ (the dual of $E$):
namely $\nabla^{E^\vee_Q}(f)=0$.

The natural evaluation map, restricted to $f$ induces a linear map
$$\eqalign{ev: H^0({X_K}, E(x))&\longrightarrow H^0(X_Q,
\O_{X_Q}(x))\cr P&\longrightarrow \langle P,f\rangle.\cr}$$ In
particular if $P\in H^0({X_K},E(x))$  and the minimal subbundle
with connection $G$ generated by $P$ has rank $r$, we will denote
by $F_i$ the sections $ev(\tilde P_i)\in H^0(X_Q;
\O_{X_Q}(x+2r))$.

The main theorem of this section is the following Zero Lemma:

\label zerolemma. theorem\par\thm Suppose that we are in the
hypotheses above, and that for every algebraic subbundle
$K\hookrightarrow E^\vee$, we have that $f\not\in H^0(X_Q, K_Q)$.
Then there exists a constant $C$ depending only on $E$, $f$ and
the fixed connections, but independent on $P$ such that
$$ord_Q(F_0)\leq x\cdot rk(G)+C.$$
\endthm

Observe that $ord_Q(F_0)=ord_Q(\langle P,f\rangle)$.

\rmk The condition on $f$ means that $f$ is not algebraically
degenerate: once one fix an algebraic trivialization of $E_F$, the
coordinates of $f$ are linearly independent over $F$.\endrmk

\Proof First of all we claim the following: $ord_Q(F_i)\geq
ord_Q(F_0)-i$: by definition
$$\eqalign{F_i&=\langle
P_i\otimes (s')^{2(r-i)};f\rangle\cr &=\langle P_i;f\rangle\otimes
(s')^{2(r-i)};\cr}$$ thus, since $s'$ do not vanish in $Q$,
$ord_Q(F_i)=ord_Q(\langle P_i;f\rangle)$. Suppose that $e$ is a
local generator of $H^{\otimes x+2i}$ and $z$ is a local
coordinate around $Q$. Then we may suppose that $\langle
P_i;f\rangle=z^a\cdot e$ for some positive integer $a$. The
evaluation map
$$ev: E(x+2i)\otimes E^\vee\longrightarrow \O(x+2i)$$
is a morphism of vector bundles with connection; thus, we may find
an analytic function $h$ in a neighborhood of Q such that
$$\eqalign{az^{a-1}he+z^a\nabla_\partial(e)&=\nabla_\partial\langle
P_i;f\rangle\cr &= \langle\nabla_\partial P_i;f\rangle+\langle
P_i;\nabla_\partial f\rangle\cr &= \langle P_{i+1};f\rangle.\cr}$$
The claim follows by induction on $i$.

We need to generalize to higher rank the notion of the order of
vanishing of a section:

Let $V$ be a vector bundle on ${X_K}$ and $f\in H^0(X_Q; V_Q)$ be
a non zero section. If we fix a trivialization of $V_Q$, we may
write $f$ as $(f_1,\dots, f_r)$ where $r$ is the rank of $V$ and
$f_i$ are power series in one variable.

\label orderofvanishing. definition\par\defi The order of
vanishing of $f$ in $Q$ is the integer $\min_i\{ ord_Q(f_i)\}$.
\enddefi

One easily sees that the order of vanishing of $f$ is independent
on the choice of the trivialization.

The theorem will be consequence of the following lemma

\label boundedfamily. lemma\par\lemma There is a constant $C$
depending only on the vector bundle with connections $E$ and $f$
with the following property: let $\cal F$ be a vector bundle with
connections and
$$\alpha: E^\vee\twoheadrightarrow\cal F$$
be a {\rm surjective} morphism of vector bundles with connections.
Let $[f]:=\alpha(f)\in H^0(X_Q,{\cal F}_Q)$; then
$$ord_Q([f])\leq C.$$
\endlemma

We first show how the lemma implies the theorem. Recall the
following standard properties of vector bundles:

a) (Cramer rule) If $G$ is a vector bundle of rank $r$ then there
is a canonical isomorphism
$$\det(G)\otimes G^\vee\simeq\bigwedge^{r-1}G.$$

b) There is a constant $C$ depending only on $E$ such that, if
$G\hookrightarrow E(x)$ is a sub bundle of rank $r$, then
$\deg(G)\leq rx+ C$.

Denote by $r$ the rank of $G$. The inclusion $G\hookrightarrow
E(x+2r)$ give rise to a surjection $\alpha:
E^\vee\twoheadrightarrow G^\vee(x+2r)$. Denote by $[f]$ the image
of $f$ in $H^0(X_Q,G^\vee(x+2r))$.

We may suppose that $\tilde P_0,\dots,\tilde P_{r-1}$ are linearly
independent elements of $G_F$ thus $\tilde P_0\wedge\tilde
P_1\wedge\dots\wedge\tilde P_{r-1}$ is a non zero global section
of $\bigwedge^r G$. Since, by property (b) above, there is a
constant $C_1$ depending only on $E$ such that
$\deg(\bigwedge^rG)\leq xr+ C$, we have that $ord_Q(\tilde
P_0\wedge\tilde P_1\wedge\dots\wedge\tilde P_{r-1})\leq xr+ C_1$.
By lemma \ref{boundedfamily} above, there is a constant $C_2$ such
that $ord_Q([f])\leq C_2$.  The isomorphism given  by  the Cramer
rule (a) give rise to the following equality:
$$(\tilde P_0\wedge\tilde P_1\wedge\dots\wedge\tilde
P_{r-1})\otimes[f]=\sum_i(-1)^i(\tilde
P_0\wedge\dots\wedge\widehat{\tilde P_i}\wedge\dots\wedge\tilde
P_{r-1})\otimes F_i;$$ thus
$$C_1+C_2+rx\geq ord_Q((\tilde P_0\wedge\tilde P_1\wedge\dots\wedge\tilde
P_{r-1})\otimes[f])\geq \inf_i\{ord_Q(F_i)\}\geq ord_Q(F_0)-r.$$ The
conclusion follows.

\rmk Observe that the constant $C$ of the theorem is sum of two
terms: the first is purely geometrical, it is essentially related
to the measure of the stability of $E$; the second term is
analytical and it is related to the structure of the specific
solution of the differential equation. \endrmk

\Proof (of Lemma \ref{boundedfamily})  We start with a
proposition:

\label bounded2. proposition\par\prop Let $V$ be a vector bundle
with singular connection on ${X_K}$; There exists a constant $C$
with the following property: Let $L$ be a line bundle with
singular connection on ${X_K}$ with a surjection $\alpha:
V\twoheadrightarrow L$ (as vector bundles with singular
connections). Then
$$\deg(L)\leq C.$$\endprop

We first show how the proposition implies the lemma: Apply the
lemma with $V=\bigwedge^r E^\vee$ and we find a constant,
depending only on $E$ such that, for every subbundle with
connection $G$ of $E$, we have that $\deg(G^\vee)\leq C$.

The degrees of the sub vector bundles with connection of $E$ are
uniformly bounded. Consequently, by the theory of the Hilbert
scheme, we can find a projective variety $\underline{Hilb_E}$, a
vector bundle $R$ on ${X_K}\times\underline{Hilb_E}$ and a
surjection $v: pr^\ast_1(E)\twoheadrightarrow R$ such that, for
every vector bundle $V$ with connection which is quotient of $E$,
there is a point $q\in\underline{Hilb_E}$ such that
$E\twoheadrightarrow V$ is the restriction of $v$ to
${X_K}\times\{ q\}$. For every $q\in\underline{Hilb_E}$ denote by
$R_q$ the vector bundle $R\vert_{{X_K}\times\{q\}}$on ${X_K}$.

Let $\underline{Hilb_E}_Q$ be the completion of
${X_K}\times\underline{Hilb_E}$ around the Cartier divisor
$\{Q\}\times\underline{Hilb_E}$. The section $f$ defines an
element of $H^0(\underline{Hilb_E}_Q, R_Q)$; thus, for every $q\in
\underline{Hilb_E}$, a global section $[f_q]$ of the localization $(R_q)_Q$ of
$R_q$ in $Q$.
Consequently we find a function
$$\eqalign{ord_Q: &\underline{Hilb_E}(\CC)\longrightarrow \Z\cr
q&\longrightarrow ord_Q([f_q]).\cr}$$ The local expression of the
function $ord_Q([f_q])$ shows that it is upper semi continuous for
the Zariski topology and since $\underline{Hilb_E}$ is compact,
the conclusion follows.

\smallskip

\Proof (of proposition \ref{bounded2}) We begin by fixing some
notation: denote by $m$ the rank of $V$. We fix a point $p$ on
${X_K}$ which is regular for the connection. Denote by $k_p$ the
completion of $\CC({X_K})$ with respect to the valuation induced
by $p$. We also fix an algebraic trivialization of $V$ near $p$.
Since the connection is regular around $p$, the space of
horizontal sections of the module with connections $V_p$ has
dimension $m$. Thus the space of {\it algebraic} horizontal
sections of $V_p^\vee$ is finite dimensional of dimension less or
equal then $m$.

Every line bundle with singular connection and quotient of $V$
defines a section $g$ of $V^\vee_p$ which is  horizontal. Thus,
$g$ belongs to a finite dimensional $\CC$--vector space $W$. The
line bundles $L$ which are quotient of $V$ are in bijection with
points of $\PP^{m-1}(\CC({X_K}))$ thus with algebraic maps
$\varphi_L:{X_K}\to\PP^{m-1}$ (modulo the action of $PGL(m)$).

Fix a basis $g_1,\dots, g_r$ of $W$. Each $g_i$ corresponds to a
quotient line bundle $L_i$ of $V$. We can associate  to every line
bundle with connection quotient of $V$, an element $g$ of $W$,
thus a linear combination of the $g_i$'s. The lemma below shows
that every line bundle quotient of $V$ which is obtained from a
linear combination of the $g_i$ has degree bounded by the maximum
of the degree of the $L_i$'s; thus the conclusion follows.

\label boundedheight. lemma\par\lemma Let
$L_i\hookrightarrow\O^m_{X_K}$ ($i=1,2$) be sub line bundles.
Consider the map
$$\eqalign{+ :\O_{X_K}^m\oplus\O_{X_K}^m &\longrightarrow \O_{X_K}^m\cr
(x,y)&\longrightarrow x+y.\cr}$$ Let $M$ be the image of
$L_1\oplus L_2$ via $+$, then $\deg(M)\geq\min\{\deg(L_i)\}$.
\endlemma

The proof of the lemma is elementary once one observe that there
is a surjection $L_1\oplus L_2\twoheadrightarrow M$.

\endsection

\

\section Nevanlinna theory and order of growth of sections\par

\

In this section we will recall the main definitions and theorems
about the order of growth of analytic maps. Most of these things
are classical, cf. for instance [GK], but the approach we have
here is a little bit different. One may may find details and
possible generalizations in [Ga].

Let $\overline{X}$ be a smooth projective curve over $\CC$ and $D$
an effective divisor on it. Let $d$ be the degree of $D$. We
denote by $X$ the affine curve $\overline{X}\setminus \{\vert
D\vert\}$.

\label growningfunction. theorem\par\thm Let $p\in X$. then, up to
an additive scalar,  there exists a unique function
$g_p:\overline{X}\to [-\infty;+\infty]$ with the following
properties:

a) it satisfies the differential equation
$$dd^cg_p=\delta_p-{{1}\over{d}}\cdot\delta_D;$$
$\delta_p$ (resp. $\delta_D$) being the dirac operator on $p$
(resp. on $D$).

b) It is a $C^{\infty}$ function on
$\overline{X}\setminus\{p\}\cup\{\vert D\vert\}$.

c) There is a open neighborhood $U$ of $p$ and an harmonic
function $v_p$ on $U$ such that
$$g_p\vert_U=\log\vert z-p\vert^2+v_p.$$
\endthm

This theorem is already proved in [Ga] in a more general
situation. We give here a sketch of proof in this case for
reader's convenience.

\Proof Fix a (Kh\"aler) metric $\omega$ on $\overline{X}$. Let
$\Delta_{\overline\partial}$ be the Laplace operator associated to
it. The operator $T:=\delta_p-{{1}\over{d}}\cdot\delta_D$ is
orthogonal to the constants. Thus there is a $(1,1)$ current
$\alpha$ on $\overline{X}$ such that
$\Delta_{\overline\partial}(\alpha)=T$. Since $T$ is smooth on
$\overline{X}\setminus\{p\}\cup\{\vert D\vert\}$, the form
$\alpha$ is also smooth on it. The operator $L:=
\cdot\wedge\omega$ induces an isomorphism between ${\cal
D}^{(0,0)}(\overline{X})$ and ${\cal D}^{(1,1)}(\overline{X})$
(${\cal D}^{(i,i)}(\overline{X})$ being the space of currents of
degree $(i,i)$). Thus there in a function $\tilde g_p$ such that
$L(\tilde g_p)=\alpha$. Since, for a suitable constant $c$, we
have that $dd^c(g)=cL(\Delta_{\overline\partial}(g))$, points (a)
and (b) are easily deduced. Point (c) is similar.

\smallskip

The functions $g_p$ are exhaustion functions in the sense of [GK]:

\label exhaustion. lemma\par\lemma For every constant $C$, we have
that $g_p^{-1}((C,+\infty])$ is a non empty neighborhood of $D$ in
$\overline{X}$.
\endlemma

\Proof Fix a metric $\Vert\cdot\Vert$ on $\O_{\overline{X}}(D)$.
Let $\II$ be the canonical section of $\O_{\overline{X}}(D)$. By
Poincar\'e--Lelong equation, the function
$g_p+{{1}\over{d}}\log\Vert\II\Vert^2$ is smooth near $D$. The
conclusion follows.

\smallskip

In the following, we will call such a function $g_p$, {\it an
exhausting function for $X$ and $p$}. Observe that an argument
similar to the one above gives

\label comparationofg. proposition\par\prop Let $p$ and $q$ points
on $X$. Let $g_p$ and $g_q$ be exhausting functions for $X$ and
$p$ and $q$ respectively. Then there is a constant $C_{p,q}$ and
an open neighborhood $V$ of $D$ such that, for every $z\in V$ we
have $$\left\vert g_p(z)-g_q(z)\right\vert\leq C_{p,q}.$$
\endprop

If $p\in X$, we fix a function $g_p$ as in the theorem above. For
every positive real number $r$, we consider the following two
closed sets of $X$
$$B(r):=\left\{ z\in X \; {\rm s.t\;}\; g_p(z)\leq\log(r)\right\}
\;\;\;\;{\rm and }\;\;\;\;S(r):=\left\{ z\in X \; {\rm s.t\;}\;
g_p(z)=\log(r)\right\}.$$

The function $g_p$ is strictly related with the Green function on
$B(r)$:

We firstly recall the definition of the {\it Green functions}:

\label greenonsurfaces. definition\par\defi Let $U$ be a regular
region on a Riemann surface $M$ and $p\in U$. A {\rm Green
function} for $U$ and $p$ is a function $g_{U;p}(z)$ on $U$ such
that:

a) $g_{U;p}(z)\vert_{\partial U}\equiv 0$ continuously;

b) $dd^cg_{U;p}=0$ on $U\setminus\{ p\}$;

c) near $p$, we have $g_{U;p}=-\log\vert z-p\vert^2+\varphi$, with
$\varphi$ continuous around $P$. \enddefi

One extend $g_{U,p}$ to all of $X$ by defining $g_{U,p}\equiv 0$
outside the closure of $U$. We easily deduce from the definitions
that $dd^cg_{U;p}+\delta_P=\mu_{\partial U;p}$ where
$\mu_{\partial U;p}$ is a positive measure of total mass one and
supported on $\partial U$.

Moreover the following is true:

\label uniquegreen. proposition\par\prop The Green function, if it
exists, it is unique.\endprop

The following gives the relation between the function $g_p$ and
the Green functions on $B(r)$:

\label greenonB. proposition\par\prop Let $r$ be a positive real
number. The function
$$g_p^r:=\log(r)-g_p\vert_{B(r)}$$
is the Green function of $B(r)$ and $p$. Consequently, for every
$p$ and $q$ in $X$ there is a constant $C$, depending on $p$ and
$q$,  such that, for every $r$ sufficiently big,
$$\left\vert g_p^r(q)-\log(r)\right\vert\leq C.$$
\endprop
The proof follows from the definitions.
\smallskip

By Stokes theorem, one easily verify that, in this case,
$\mu_{S(r);p}$ is the positive measure $d^cg_p\vert_{S(r)}$.

Let $Z$ be a projective variety and $L$ be an ample line bundle on
it equipped with a positive metric. Denote by $c_1(L)$ the first
Chern form associated to it.

Let $\gamma :X\to Z$ be an analytic map. We define the height
function associated to it:
$$T_\gamma(r):=\int_0^r{{dt}\over{t}}\int_{B(t)}\gamma^\ast(c_1(L))=\int_Xg_p^r\cdot\gamma^\ast(c_1).$$
The order of growth of the map $\gamma$ is given, as
$$\limsup_{r\to+\infty} {\log{T_\gamma(r)}\over \log(r)}.$$

More generally, if $M$ is an hermitian line bundle on $X$, we
define
$$(M,X)(r):=\int_0^r{{dt}\over{t}}\int_{B(t)}c_1(M)=\int_Xg_p^r\cdot c_1(M).$$

Some remarks are necessary, we can find the proofs for instance in
[Ga]:

-- The order of growth is independent on:

(i) The choice of the ample line bundle $L$ and on the metric on
it.

(ii) The choice of the point $p$.

-- If $\gamma$ is the inclusion in $\overline{X}$, or more
generally if $\gamma$ is an algebraic map (cf. [GK]), then, there
is a constant $C$ such that \labelf rationalfunc\par$$\left\vert
{{T_{\gamma}(r)}\over{\log(r)}}\right\vert\leq C.\eqno{(\numfo)}$$

-- The Stokes and Poincar\'e--Lelong formulas give rise to the
first main theorem: Let $Y\in H^0(X,M)$ be a global section. We
define {\it counting function} of $Y$: suppose that $div(Y)=\sum
n_z z$ (the sum may possibly be infinite), and to simplify, that
$p\not\in div(Y)$, then
$$N_Y(r):=\int_Xg_p^r\cdot\delta_{div(Y)}=\sum_{g^r_p(z)<\log(r)}n_zg_p^r(z).$$
The {\it First Main Theorem} (FMT) holds:
$$N_Y(r)-\int_{S(r)}\log\Vert
Y\Vert^2\mu_{S(r),p}=(M,X)(r)+\log\Vert Y\Vert^2(p).$$

The term $-\int_{S(r)}\log\Vert Y\Vert^2\mu_{S(r),p}$ is often
denoted by $m_Y(r)$ and called {\it proximity function of $Y$}.

\smallskip

Let $E\to\overline{X}$ be an hermitian vector bundle and
$p:\PP:=\underline{Proj}(\O\oplus E^\vee)\to \overline{X}$ be the
associated compactification of it. Let $\MM$ be the tautological
line bundle of $\PP$; since $E$ is hermitian, $\MM$ is naturally
equipped with the relative Fubini--Study metric. The surjection
$\O\oplus E^\vee\to E^\vee$ defines an inclusion
$\PP(E)\hookrightarrow\PP$ (the divisor at infinity) and the image
is a global section of $\MM$. It is well known that, if $M$ is a
sufficiently ample line bundle on $\overline{X}$ then $\MM\otimes
p^\ast(M)$ is a very ample line bundle on $\PP$.

Let $f:X\to E$ be an analytic section of $E$. It canonically
defines an analytic map $f_\PP:X\to \PP$. By definition, the order
of growth of $f_\PP$ is
$\limsup_{r\to\infty}{{\log(f_\PP^{\ast}(\MM\otimes
p^\ast(M));X)(r)}\over{\log(r)}}$. Observe that by
\ref{rationalfunc} the order of growth of $f_\PP$ is independent
on $M$.

\label orderofgrowth. definition\par\defi We define the order of
growth of the section $f$ to be the number
$$\rho:=\limsup_{r\to\infty}{{\log(f_{\PP}^{\ast}(\MM);X)(r)}\over{\log(r)}}.$$

\enddefi

\label order1. lemma\par\lemma Suppose that $f$ is a section of
order strictly less then $\rho$. Then there is a constant $C$ such
that
$$\int_{S(r)}\log\Vert f\Vert\mu_{S(r),p}\leq Cr^\rho.$$
\endlemma

\Proof Observe that $f_\PP$ do not intersect $\PP(E)$ and that, if
$q\not\in\PP(E)$, then $\Vert \PP(E)\Vert^2(q)={{1}\over{1+\Vert
q\Vert^2}}$. Thus, by FMT, there is a constant $C$ such that
$$(f_\PP^{\ast}(\MM);X)(r)={{1}\over{2}}\int_{S(r)}\log\left(1+\Vert
f\Vert^2\right)\mu_{S(r),p}+ C.$$ The conclusion easily follows.

\smallskip

We will show that, given a section with finite order of growth,
and two points, we can estimate the size of a related section in
one point, if we know that this vanishes to an high order on the
other point.

Fix two points $p_1$ and $p_2$ in $X$.

Suppose that $E$ is an algebraic vector bundle over $\UX$. Fix an
ample line bundle $H$ on $\UX$. We suppose that $E$ and $H$ are
equipped with smooth metrics. For every positive integer $x$,
denote by $E(x)$ the vector bundle $E\otimes H^{\otimes x}$.

Fix an {\it analytic section} $f\in H^0(X, E^\vee)$ having order
of growth $\rho$. For every $P\in H^0(\UX, E(x))$ denote by $F$
the analytic section section $f(P)\in H^0(X, H^{\otimes x})$.

We will show that one can bound the size of $F$ in $p_2$ in terms
of the sup norm of $P$, the order of vanishing of $F$ in $p_1$ and
the order of growth of $f$.

\label estimateinp2. theorem\par\thm There is a constant $c_1$
depending only on $H$, a constant $c_2$ depending only on $f$ and
a constant $c_3$ depending only on $p_1$ and $p_2$, for which the
following holds:

For every section $P\in H^0(\UX; E(x))$ such that

-- $\log\sup\{\Vert P\Vert\}\leq B$;

-- $ord_{p_1}(F)\geq Ax-b$;

We have the following estimate, for every $x\gg 0$:
$$\log\Vert F\Vert(p_2)\leq B-{{Ax}\over{\rho}}\cdot\log(x)+
c_2x+{{c_1}\over{\rho}}\cdot x\log(x).$$ \endthm

Observe that the constant $c_1$ depends only on $H$.

\Proof First of all remark that, by Cauchy--Schwartz inequality,
we have that $\Vert F\Vert\leq\Vert P\Vert\cdot\Vert f\Vert$. By
Stokes formula we have that, for every real number $r$,
$$\int_X\log\Vert F\Vert\cdot dd^cgg_{P_2}^r=\int_Xdd^c\log\Vert
F\Vert\cdot g_p^r.$$ By the definition of Green function, the Left
Hand Side is
$$\int_{S(r)}\log\Vert F\Vert\mu_{S(r), p}-\log\Vert F\Vert(p_2).$$
By the Cauchy--Schwartz inequality, the hypotheses and Lemma
\ref{order1} the first term of the sum is bounded by
$$B+c_2r^\rho.$$
Since $ord_{p_1}(F)\geq Ax$, by \ref{greenonB}, and the
Poincar\'e--Lelong formula, the Right Hand Side is surely bigger
then
$$(Ax-b)\log(r)-x(H,X)(r);$$
Since $H$ is algebraic, with a metric smooth at infinity, the last
term of this sum is surely lower bounded by $-x\cdot c_1\log(r)$,
for a suitable $c_1$ depending only on $H$. The conclusion follows
by taking $r=x^{1/\rho}$.

\endsection

\

\section Order of growth at finite places\par

\

In this section we will recall the definitions and the principal
properties of the $LG$ germs. This definition is given in [Ga] and
there developed in a greater generality, and here we just recall
it (and explain in the special situation we need it) for reader
convenience. The notion of $LG$ germ is similar to the notion of
$E$ function developed by Siegel, Shidlowski and others. When we
are in presence of $LG$ germs, we can estimate the order of growth
of sections at all the finite places at the same time. It is our
opinion that, the notion of $LG$ germs and the order of growth of
sections (or more generally of analytic maps) are two concepts
which may be in contrast; and from this contrast we may deduce non
trivial results.

Let $K$ be a number field, $O_K$ its ring of integers and $M_K$
the set of places of $K$. We will denote by $M_{fin}$ the set of
finite places of $K$. If $v\in S$,v we denote by $K_v$ the
completion of $K$ with respect to $v$ and by $O_v$ its ring of
integers; if $M$ is an $O_K$ module, we denote by $M_v$ the $K_v$
vector space $M\otimes K_v$ and by $M_{O_v}$ the $O_v$ module
$M\otimes O_v$.

Let $X_K$ be a smooth projective curve over $K$, $H_K$ be an ample
line bundle over it. Let $E_K$ be a vector bundle of rank $m$ over
$X_K$. We will use the same notation then before. Denote by
$E^\vee_K$ the dual of $E_K$

Let $\X\to\Spec(O_K)$ be a model of $X_K$. We can (and we will)
suppose that $H_K$ (resp. $E_K$) extends to a line bundle $H$
(resp. to a vector bundle $E$) over $\X$.

Let $p_K\in X_K(K)$ be a rational point. It extends to a section
$p:\Spec(O_K)\to \X$. denote by $\hX_p$ the completion of $X_K$
near $p_K$. We denote $\hat{\X}_p$ the completion of $\X$ near
$p$. Denote by $H_p$, $E_p$ etc (resp $H_{p,K}$, $E_{p;K}$ etc.)
the restriction of $H$, $E$ etc to $\hat{\X}_p$ (resp. of $H_K$,
$E_K$ etc to $\hX_p$). Up to pass to a finite extension of $K$, we
may suppose that $\hat{\X}_p$ is isomorphic to
$Spf(O_K[\![Z]\!])$; we fix such an isomorphism; we may also
suppose that $H_p$, and $E_p$ are isomorphic to the trivial bundle
of the corresponding rank; we fix such isomorphisms.

Let $f\in H^0(\hX_p;E_{p,K}^\vee)$.

Since we fixed the isomorphism of $E^\vee_p$ with
$\O_{\hat{\X}_p}^m$, the section $f$ can be written as $m$ power
series. Essentially, $f$ will be a $LG$ germ, if we can control
the denominators of these power series in terms of powers of the
factorials.

Once we fixed the isomorphisms above, the section $f$ may be
written as $m$ power series
$$f=\left(\sum_{i=1}^{\infty}a_i(1)Z^i;\dots
;\sum_{i=1}^{\infty}a_i(m)Z^i\right);$$ with $a_i(j)\in K$.

\label lggerms. definition\par\defi We will say that $f$ is a
$LG$--germ of type $\alpha$ if the following holds:

a) For every place $v\in M_K$ and every $j=1,\dots, m$, the power
series $\sum_ia_i(j)Z^i$ have positive radius of convergence;

b) There is a finite set of places $S$ such that, if $v\not\in S$
there is a constant $C_v$ such that, for every $j=1,\dots, m$,
$$\Vert a_i(j)\Vert_v\leq {{C^i_v}\over{\Vert i!\Vert_v^\alpha}};$$

c) $\prod_{v\not\in S}C_v <\infty.$
\enddefi

Following the proofs of [Ga] \S 3, one may prove that:

-- The notion of $LG$--germ of type $\alpha$ do not depend on the
chosen isomorphisms; thus the notion depends only on the germ of
section. However, notice  that the constants $C_v$'s may depend on
the choices.

-- If $E$ is equipped with a connection which is regular at $p$; a
formal horizontal section is a $LG$--germ of type 1.

-- If moreover, for almost all $v\in M_K$, the connection has
vanishing $p$--curvature, then the formal horizontal section is an
$LG$--germ of type zero.

The last two sentences are essentially proved in [Bo].

\smallskip

Suppose we fixed $s$ points $p_1,\dots, p_s$ in $X_K(K)$. Suppose that,  for
every point $p_j$ we have an $LG$ germ $f_j\in H^0(\hX_{p_j};E_{p_j,K}^\vee)$ of
type $\alpha$. If we take a suitable blow up of $\X$, we may suppose that the
$p_j$'s extend to sections $P_j:\Spec(O_K)\to \X$ which do not intersect.

For every integer $x$, denote by $G_x$ the $O_K$--module $H^0(\X;
E\otimes H^{\otimes x})$. For every $j$, the  The section $f_j$ induces a
$O_K$--linear
map $\langle\cdot;f_j\rangle: G_x\to H^0(\hX_{p_j};H_{p_j,K}^{\otimes
x})$; and, by composition a map
$$\langle\cdot;f\rangle:=(\langle\cdot;f_1\rangle,\langle\cdot;f_2\rangle,
\dots ,
 \langle\cdot;f_s\rangle):G_x\longrightarrow
\oplus_{j=1}^rH^0(\hX_{p_j};H_{p_j,K}^{\otimes
x})$$

Let $j\in\{1,\dots, s\}$. For every positive integer $i$, denote by
$\hat{\X}_{p_j}^i$ (resp.
$\hX^i_{p_j}$) the $i$-th infinitesimal neighborhood of $p_j$ (resp.
$p_{j,K}$) in $\X$ (resp. $X_K$). Similarly we denote by $H_{p_j,i}$
etc. the restriction of $H$ etc. to $\hat{\X}_{p_j}^i$. Since $\X$ is
smooth near each of $p_j$'s, the sheaf of Kh\"aler differentials
$\Omega^1_{\X/O_K}$ is locally free in a neighborhood of $p_j$.
Denote by $T_{p_j}\X$ the restriction to $p_j$ of the dual of it.

Denote by $res_j:\oplus_{j}H^0(\hX_{p_j};H_{p_j,K}^{\otimes x})\to
\oplus_{j}H^0(\hX^i_{p_j};H_{p_j,i}^{\otimes
x})\otimes K$ the restriction map, by $\langle\cdot;f\rangle_i$ the  map
obtained by composing
$\langle\cdot;f\rangle$ with the $res_j$ and by
$G_x^i$  the kernel of $\langle\cdot;f\rangle_i$.

The exact sequence
$$0\to \oplus_{j=1}^sH^0(p_j, H^{\otimes x}\otimes (T_{p_j}\X)^{\otimes
-i})\to
\oplus_{j=1}^sH^0(\hX_{p_j}^{i+1};H^{\otimes x}_{p_j,i+1})\to
\oplus_{j=1}^sH^0(\hX_{p_j}^i;H^{\otimes x}_{p_j,i})$$ and the snake lemma,
induces a canonical
inclusion
$$\gamma_x^i: G_x^i/G_x^{i+1}\longrightarrow \oplus_{j=1}^sH^0(p_j, H^{\otimes
x}\otimes (T_{p_j}\X)^{\otimes
-i})\otimes K.$$

For every finite place $v\in M_{fin}$ both $(G_x^i/G_x^{i+1})_v$
and $\oplus_{j}H^0(p_j, H^{\otimes x}\otimes (T_{p_j}\X)^{\otimes -i})_v$ are
equipped with norms, induced by the integral structure. Observe that,
naturally the former has the $\sup$ norm. Thus we
may compute the norm $\Vert \gamma_x^i\Vert_v$ of the operator
$\gamma_x^i$.

When $f$ is an $LG$--germ, then one can bound the norms at the
finite places of the $\gamma_x^i$'s.

\label boundofnormsoflggerms. theorem\par\thm With the notations
as above, suppose that $f$ is a $LG$--germ of type $\alpha$ then
there exists a constant $C$ such that
$$\sum_{v\in M_{fin}}\log\Vert\gamma_x^i\Vert_v\leq[K:{\Bbb
Q}]\alpha\cdot i\cdot\log(i)+C(i+x).$$
\endthm

\Proof We first remark the following general statement: Let $k$ be
a normed field and $\varphi: V^1_k\to V^2_k$ be a linear map
between finite dimensional normed vector spaces over $k$. Let
$V^i\subset V_k^i$ be the set of elements of norm less or equal
then one. Suppose that there exists a constant $A\in k^\ast$ such
that, for every element $v\in V^1$ we have that
$\varphi(Av)\subset V^2$, then
$\Vert\varphi\Vert\leq{{1}\over{\Vert A\Vert}}$.

For every $j\in\{1,\dots, s\}$, let $pr_j:\oplus_{j}H^0(p_j, H^{\otimes
x}\otimes (T_{p_j}\X)^{\otimes -i})\to H^0(p_j, H^{\otimes
x}\otimes (T_{p_j}\X)^{\otimes -i})$ be the projection. Fix such a $p_j$.
Let $v\not\in S$. The restriction of $\hat\X_{p_j}$ to $\Spec(O_v)$ is
isomorphic to $Spf(O_v[\![Z]\!])$ (via an  isomorphism fixed as 
above).

Suppose that $P\in (G_x^i)_{O_v}$. Since we fixed an isomorphism
of  $E_{p_j}$ with the trivial vector bundle of rank $m$, the
restriction of $P$ to $(\hat\X_{p_j})_v$ is represented by
$(g_i,\dots, g_m)$ with $g_i\in O_v[\![Z]\!]$. By definition
$\langle P,f_j\rangle =\sum_{s=1}^mg_s\sum_\ell
a_\ell(s)Z^\ell=h_j(Z)$. Since $P\in (G_x^i)_{O_v}$, we have that
$h_j(Z)=\sum_{\ell=i}^\infty h_\ell Z^\ell$ and $pr_j\circ\gamma_x^i(P)=h_i$.
Since $f_j$ is a $LG$--germ of type $\alpha$, we have that
$${{\Vert i!\Vert_v^\alpha}\over{C_v^i}}\cdot \Vert h_i\Vert_v\leq
1.$$ 

The norm on $\oplus_{j}H^0(p_j, H^{\otimes
x}\otimes (T_{p_j}\X)^{\otimes -i})$ is the $\sup$ of the norms on each
factors, consequently, for every $v\not\in S$ we have
$${{\Vert i!\Vert_v^\alpha}\over{C_v^i}}\cdot\Vert \gamma_x^i(P)\Vert_v\leq 1.$$

The conclusion follows from the remark at the beginning of this proof, 
the Stirling formula and the
standard Cauchy inequality at the places in $S$.

\endsection

\

\label Esections. section\par\section $E$--sections of type
$\alpha$\par

\

In this section we will introduce the concept of $E$-- sections of
type $\alpha$ of a vector bundle over an affine curve. These are
analytic sections of an algebraic vector bundle whose order of
growth is the inverse of the type of their formal development in a
fixed algebraic point. The main examples of $E$--sections are the
$E$--functions of the theory of Siegel--Shidlowski (they are of
type 1) or the more recent "arithmetic series of order $s$"
introduced by Andr\'e in [An1].

Let $\UX_K$ be a smooth projective curve defined over a number
field $K$. Fix $s$ points  $p_1,\dots, p_s\in \UX(K)$.

Let $E_K$ be a vector bundle over $\UX$ of rank $m$.

\label esection. definition\par\defi A section $f:=(f_1,\dots,f_s)\in
E_{p_1,K}\oplus E_{p_2,K}\oplus\dots\oplus E_{p_s,K}$ is
said to be a {\rm $E$--section of type $\alpha$} if the following
holds:

a) For every $j$,  the germ of section $f_i\in E_{p_j,K}$ is an $LG$--germ of
type $\alpha$.

b) There exists a non empty subset $S_K\subseteq M_{\infty}$ of
cardinality $a$ such that the following holds:  for every
$\sigma\in S$ there is an affine open subset $U_\sigma$ of
$\UX_\sigma$ and an {\rm analytic} section $f_\sigma\in
H^0(U_\sigma; E_\sigma)$ such that :

\noindent (b.1) For every $j$, the germ of $f_\sigma$ at $p_j$ is $f_j$;

\noindent(b.2) the section $f_\sigma$ has order of growth
$\rho_\sigma={{a\cdot s}\over{[K:\QQ]\alpha}}$.

\enddefi

\rmk (a) An $E$--function in the sense of Siegel Shidlowski is an
$E$-- section of type 1; in this case we have only one point and $S_K=
M_\infty$;

(b) An "arithmetic Gevrey series of order $s<0$" in the sense of
Andr\'e [An1] is a $E$--section of type $-s$, again we only have one point
and
$S_K=M_\infty$.

(c) If $L/K$ is a finite extension and $f\in E_{p,K}$ is a $E$
section of type $\alpha$ over $K$, then $f\in E_{p,L}$ is a
$E$--section of type $\alpha$; take as $S_L$ the set $\tau$ such
that $\tau/\sigma$ for $\sigma\in S_K$.

(d) Notice that, on the projective line, the main differences
between $E$--functions and $E$--sections are: (1) $E$ sections may
have order of growth which is not one. (2) (more important)
$E$--sections may have finitely many essential singularities,
whereas $E$--functions are always entire functions.

(e) An interesting example of $E$ section, and our main theorem will concern
such example, is given by an horizontal section of a fibre bundle with
meromorphic connections having order of growth $\rho$ and assuming 
algebraic  values at $\rho$ regular rational points: in each of the rational
points the section will be an $LG$--germ of type 1.  

(f) By Corollary \ref{corollaofsl} the order of growth cannot be less then $s$.
\endrmk

In this section we show that, given an $E$--section, it is
possible to construct sections with high order of vanishing and
bounded sup norm.

First of all we have to fix integral structures: As in the
previous section, we suppose that $\X\to\Spec(O_K)$ is a regular
projective model of $\UX_K$ and $E_K$ extends to a vector bundle
over $\X$. We suppose also that $H$ is a relatively ample line
bundle on $\X$. For every place $\sigma\in M_\infty$, we suppose
that $E_{\sigma}$ and $H_\sigma$ are equipped with smooth metrics
(and the metric on $H$ is sufficiently positive). We also fix
metrics on $(\UX)_\sigma$. Thus, for every integer $x$, the vector
bundle $E^\vee(x)$ is an hermitian vector bundle over $\X$.

For every integer $x$, the $O_K$--module $H^0(\X, E^\vee(x))$ is
equipped with a structure of {\it hermitian $O_K$--module}: for
every $\sigma\in M_{\infty}$, $H^0(\UX_\sigma; E^\vee(x)_\sigma)$
is equipped with the $L^2$ metric (notice that the $L^2$ norm and
the $\sup$ norm are comparable by for instance [Bo] \S 4.1). As in
the previous section, in the sequel, we will denote it by $G_x$.

Suppose we fixed an $E$--section $f\in \oplus_j E_{p_j,K}$ of type $\alpha$.
Using the notations of the previous section, for every positive
integer $i$, we obtain a natural $O_K$--linear map $G_x\to
\oplus_jH^0(\hX^i_{p_j};H_{p_j,i}^{\otimes x})\otimes K$. Again denote by
$G_x^i$ its kernel.  Put $c:=\deg(H_K)$. We want to prove that,
under these conditions, for every $\epsilon\in(0,1)$ there is a
non vanishing section of bounded norm in $G_x^{x{{c}\over{s}}m(1-\epsilon)}$.

\label smallsections. theorem\par\thm Suppose that we fixed the
hypotheses as above. Fix $\epsilon\in (0,1)$. Then, we can find a
constants $a$, depending only on $\epsilon$, $\alpha$ the points $p_j$, but
independent on the vector bundle $E$ and a constant $b$ for which
the following holds: for every sufficiently big positive integer
$x$ there is a {\rm non zero} section $P\in G_x^{x{{c}\over{s}}m(1-\epsilon)}$
such that
$$\sup_{\sigma\in M_{\infty}}\{\log\Vert P\Vert_\sigma\}\leq a\cdot x\cdot\log(x)+b\cdot
x.$$
\endthm

Before we start the proof, we need to recall some classical tools
form Arakelov geometry.

-- Suppose that $E_1$ is an hermitian $O_K$--module and $L_1,\dots, L_s$ are
hermitian line bundles over $\O_K$. Let $\varphi:E_1\to \oplus_jL_j\otimes K$ is
an injective linear map. For
every place $v\in M_K$ (finite or infinite) we denote by $\Vert
\varphi\Vert_v$ the norm of $\varphi$. One easily find that, if
$\varphi$ is non zero, then
$$\widehat{\deg}(E_1)\leq
rk(E_1)\cdot\left(\sup_j\{\widehat{\deg}(L_j)\}+\sum_{v\in
M_K}\log\Vert \varphi\Vert_v\right).$$

-- There exists a constant $\chi(K)$ depending only on $K$ such
that the following holds:  Suppose that $E$ is an hermitian
$O_K$--module. Suppose that $\widehat{\deg}(\overline E)\geq A$
then there exists a {\it non zero} element $x\in E$ such that
$$\sup_{\sigma\in M_{\infty}}\{\Vert x\Vert_{\sigma}\}\leq
-{{A}\over{rk(E)}}+\log(rk(E))+\chi(K).$$ (cf. [BGS] thm. 5.2.4
and below it).

-- if $x$ is sufficiently big, we may suppose that
$\widehat{\deg}(G_x)\geq 0$.

We can now start the proof of the theorem.

\Proof As in the previous section, for every integer $i$,  we have
an injective map 

\advance\ssnu by-1\labelf
filtration\par$$\gamma_x^i:G^i_x/G_x^{i+1}\longrightarrow
\oplus_{j=1}^rH^0(p,H^{\otimes x}\otimes(T_{p_j}\X)^{\otimes -i})\otimes K
.\eqno{{(\numfo)}}$$\advance\ssnu by1  Remark that all the
$G_x^i$'s are hermitian $O_K$ modules.  From the properties listed
above we find that we can find a constant $A$ depending only on
$p$ and $H$ such that
$$\widehat{\deg}(G_x^{i+1})\geq \widehat{\deg}(G_x^i)-rk(G_x^i/G_x^{i+1})\left(
A(x+i)+\sum_{v\in M_K}\log\Vert \gamma_x^i\Vert_v\right).$$ Since
$f$ is an $E$--germ of type $\alpha$, by theorem \ref{boundofnormsoflggerms},
there is a constant $c$ such
that
$$\sum_{v\in M_{fin}}\log\Vert\gamma_x^i\Vert_v\leq
[K:\QQ]\alpha\cdot i\cdot\log(i)+c\cdot (i+x).$$ Let $S_K$ be the
set of infinite places involved in the definition of $E$--germ. By
the classical Cauchy inequality there is a constant $C$ such that,
if $\sigma$ is an infinite place not contained in $S_K$, then
$$\log\Vert\gamma_x^i\Vert_\sigma\leq C(x+i).$$
In order to estimate the norm in the places of $S_K$ we need a
refinement of Theorem \ref{estimateinp2}.

Let $\sigma\in S_K$. Let $j\in\{ 1,\dots, s\}$, we put on the line bundle
$\O_{U_\sigma}(p_j)$
the following metric: Suppose that $\II_{p_j}$ is the canonical
section of $\O_{U_\sigma}(p_j)$, then we define $\Vert
\II_{p_j}\Vert(z)=\exp({{1}\over{2}}g_{p_j}(z))$. By adjunction, this
defines a norm on $T_{p_j}\UX_\sigma$. Let $s\in (G_i)_\sigma$; then
$f_\sigma(s)\cdot\prod_j\II_{p_j}^{-i}$ is an holomorphic section $\tilde F$ of
$(H^{\otimes x}(-\sum_jip_j))_\sigma$. To compute the norm of $\gamma_x^i$ at
the places in $S_K$ we have to  compare $\Vert\tilde
F\Vert(p_j)$, for every $j$,  with $\Vert s\Vert_{\infty}$.

By Stokes formula we find, for every real number $r$,
$$\int_{U_\sigma}\log\Vert\tilde F\Vert\cdot dd^c
g_p^r=\int_{U_\sigma}dd^c\log\Vert\tilde F\Vert\cdot g_p^r.$$

By proposition \ref{comparationofg} as soon as $r\gg 0$ we may suppose that, if
$g_{p_j}(z)\geq r$, then $g_{p_{j_i}}(z)=g_{p_{j}}(z)(1+\epsilon(z))$ with
$\vert\epsilon(z)\vert\leq\epsilon$
Thus, by the property of the Green function recalled in \S 3, the
definition of the norm on $\O(p_j)$, the Cauchy--Schwarz inequality
and the Poincar\'e--Lelong formula we find
$$\log\Vert s\Vert_{\infty}+\int_{S(r)}\log\Vert
f_\sigma\Vert d^cg_p-is(1-\epsilon)\log(r)-\log\Vert\tilde F\Vert(p)\geq
-x(H;U_\sigma)(r).$$ Thus, we can find constant $C$ and $\epsilon>0$, depending
only
on $H$ and a constant $\lambda_\sigma$ depending on $f$, such that, as soon as
$r\gg0$
$$\log\Vert\gamma_x^i\Vert_\sigma\leq -is(1-\epsilon)\log(r)+x\cdot
C\log(r)+\lambda_\sigma r^{{{as}\over{[K:\QQ]\alpha}}(1+\epsilon)}.$$

For each $i$ we put $r=i^{{\alpha[K:\QQ]}\over{a}s(1+\epsilon)}$ and we obtain
that there are constants $C_1$ $C_2$ and $\epsilon_1$ with $C_1$
depending only on $\alpha$
and $H$,  and $\epsilon_1$ as small as we want, in particular independend on
$E$,   $i$ and $x$ such that
$$\sum_{v\in M_K}\log\Vert \gamma_x^i\Vert_v\leq
xC_1\log(i)+ \epsilon_1 i\log(i)+C_2(i+x).$$

Observe that $rk(G_x^{i}/G_x^{i+1})\leq s$ consequently we can find
constants $C_i$  with $C_4$ depending only on $H$, $s$  and $\alpha$
such that, summing up all together we obtain
$$\widehat{\deg}(G_x^{x{{c}\over{s}}m(1-\epsilon)})\geq C_3x^2 -
C_4\left(\sum_{i=1}^{x{{c}\over{s}}m(1-\epsilon)}
x\log(i)+\epsilon_1i\log(i)+C_2(i+x)\right).$$
Thus, since we can take $\epsilon_1$ very small compared to $m^2$,  there are
constants $C_6$ and $C_7$ with $C_6$ independent on $E$ and $C_7$ depending on
$m$, $f$ and $H$ such
that
$$\widehat{\deg}(G_x^{x{{c}\over{s}}m(1-\epsilon)})\geq -
\left(C_6m\cdot x^2\cdot\log(x)
+ C_7\cdot x^2\right).$$

By Riemann--Roch theorem, $rk(G_x)$ is about $cm x$. By the
filtration \ref{filtration}, the rank of $G_x^{x{{c}\over{s}}m(1-\epsilon)}$ is
bigger then $\epsilon mx$. Consequently, there is a {\it non zero
section} $P$ of $G_x^{x{{c}\over{s}}m(1-\epsilon)}$ such that
$$\sup_{\sigma\in M_\infty}\{\log\Vert P\Vert_{\sigma}\}\leq
{{mC_8}\over{m\epsilon}}\cdot x\cdot\log(x)+
C\cdot x.$$ The conclusion follows since, $C_8$ depend only on the $p_j$,
$H$ and it is independent on $E$.
\smallskip

We suppose that we fixed
$s$ rational points and a fibre bundle $(E,\nabla)$ with a meromorphic
connection. We suppose that we fixed an $E$ section of type $\alpha$ in the
in the neighborhood of these points which is horizontal for the connection.
Using the connection, we may now take derivatives of $P$ in order to construct
other sections with the same properties. 

First of all we have to assure that the derivative of an integral
section is again an integral section: we may extend the connection
$\nabla :E^\vee(x)_K\to E^\vee(x)_K\otimes\Omega^1_{\hX}(D_K)$ to
an {\it integral connection}:
$$\nabla: E(x)\to E\otimes\omega_{\X/O_K}(D+V)$$
where $V$ is a vertical divisor. We also fix an integral element
$\partial\in T_{\X/O_K}(D)$ which, generically, do
not vanish at the points $p_i$'s, in case we may suppose the degree of $D$ very
big. By construction, if $P\in H^0(\X, E^\vee(x))$, then
$\nabla_\partial(P)\in H^0(\X,E^\vee(x+2+V))$; in particular
notice that $\nabla_\partial(P)$ is an integral element.For every point $p_j$,
the order
of vanishing of $\langle \nabla_\partial(P),f_j\rangle$in $p$ is one
less then the order of vanishing in $p$ of $\langle P,f\rangle$ in
$p$. Moreover a straightforward application of the classical
Cauchy--Schwarz inequality implies that, for every complex
embedding $\sigma$, the linear map $(\nabla_\partial)_\sigma
:E^\vee(x)_\sigma\to (E^\vee_\sigma(x+2))_\sigma$ has bounded
norm. Thus we proved:

\label globalsections. proposition\par\prop There is a constant
$A$ depending only on $(E, \nabla)$ and $\partial$ such that the
following holds: if $P\in H^0(\X, E^\vee(x))$ is an integral
section such that $\sup_{\sigma}\{\log\Vert P\Vert_{\sigma}\}\leq
C$ and, for every $j$ we have $ord_{p_j}(\langle P,f_j\rangle)\geq C_1$ then:

(i) $\nabla_\partial(P)$ is an integral section of $E^\vee(x+2+V)$
such that
$$\sup_{\sigma}\{\log\Vert \nabla_\partial(P)\Vert_{\sigma}\}\leq C+A;$$

(ii) For every $j$ we have $ord_{p_j}(\langle \nabla_\partial(P),f_j\rangle)\geq
C_1-1.$
\endprop

\endsection

\section Proof of the main theorem\par

\

In this section we will show how to generalize the Siegel
Shidlovski theory to an arbitrary curve and to a connection with
arbitrary meromorphic singularities and $E$ sections over an arbitrary set of
points.

\label SSoncurves. theorem\par\thm Let $\UX/_K$ be a smooth
projective curve defined over the number field $K$. Let $D_K$ be
an effective divisor on $\UX$ and $(E_K,\nabla_K)$ be a fibre
bundle of rank $m>1$ with connection with meromorphic poles along
$D_K$. Let $p_1,\dots, p_s\in \UX(K)$ be rational points.

Let $f:=(f_1,\dots, f_s)\in \oplus_j E_{p_j,K}$ be an {\rm horizontal section}
which is an
$E$--section of type $\alpha$.

Let $S_K$ be the subset of infinite places of $K$ involved in the
definition of $f$ and $\sigma\in S_K$.

Suppose that $q\in\UX\setminus\{ D,p_1,	\dots, p_s\}$. Then
$$Trdeg_K(K(f_{\sigma}(q)))=m.$$
\endthm

\rmk If we apply the theorem to $\PP^1$ with $s=1$  and we suppose that the
horizontal section is an $E$ function, we find the classical
theorem of Siegel and Shidlovski cf. [La]. \endrmk

\Proof We fix models $\X$ of $\UX$, $D$ of $D_K$ and $(E,\nabla)$ of
$(E_K,\nabla_K)$. We also fix a positive metric on the ample line
$H:=\O(D) $. Let $c$ be the degree of $H_K$, adding some points to $D$ if
necessary we may suppose that $c>s$. We eventually fix an
integral derivation $\partial\in H^0(\X,
(\omega^1_{\X/O_K})^\vee(D))$ which do not vanish at the points $p_j$'s and
$q$; notice that this can be done once we suppose that $c$ is big compared to
$s$.

We want to apply \ref{linearindependencecriterion2}, thus we need
$m$ linearly independent sections of $E^\vee(x)_q$ with satisfying
the hypotheses of loc. cit.

By Theorem \ref{smallsections}, for every $x\gg 0$ we may
construct $P_1\in H^0(\X; E^\vee(x))$ such that, using the
notations of loc. cit., for every $j$ we have  $ord_{p_j}(F_1)\geq
x{{c}\over{s}}m(1-\epsilon)$ and $\sup_
{\sigma}\{\log\Vert P_1\Vert_{\sigma}\}\leq ax\log(x)+c_1x$.

Let $P_i:=\nabla_\partial P_{i-1}$.  Applying proposition
\ref{globalsections}, we construct then $m$ integral sections
$P_1,\dots ,P_m$ such that (again using notations of loc. cit.)
$\sup_ {\sigma}\{\log\Vert P_i\Vert_{\sigma}\}\leq ax\log(x)+c_2x$
and $ord_p(F_i)\geq x{{c}\over{s}}m(1-\epsilon)-m$.

Since we suppose that $c>s$ and $x\gg 0$, we may apply The Zero Lemma
\ref{zerolemma} and we obtain that $P_1, \dots, P_m$ are linearly
independent over $K(\UX)$. Observe that $P_i\in H^0(\X,
E^\vee(x+2m+ V))$ for some fixed vertical divisor $V$.
Consequently there is a constant $c_3$ such that
$$\deg(P_1\wedge\dots\wedge P_m)=mcx+c_3.$$
By Cramer rule,
$$(P_1\wedge\dots\wedge P_m)\otimes
f=\sum_i(-1)^i(P_1\wedge\dots\wedge\hat{P_i}\wedge\dots\wedge
P_m)\otimes F_i;$$ thus, for every $j$,  $ord_{p_j}(P_1\wedge\dots\wedge
P_m)\geq x{{c}\over{s}}m(1-\epsilon)-m$; consequently, constants $\epsilon_1$
and  $c_4$ such
that
$$ord_q(P_1\wedge\dots\wedge P_m)\leq cx\epsilon_1+c_4.$$

\label computationatq. lemma\par\lemma With the notations as
above, there are constants $c_i$'s independent on $E$ (in
particular independent on $m$) and constants $b_j$ for which the
following holds: for every $x\gg 0$ there exist $m$ indices
$\ell_1<\ell_2<\dots <\ell_m$ with $\ell_m\leq m+c_1x+c_2$ such
that:

(a) $P_{\ell_1}\wedge\dots\wedge P_{\ell_m}$ is an integral
section {\rm non vanishing at $q$};

(b) $\sup_{\sigma}\{\log\Vert P_{\ell_i}\Vert_{\sigma}\}\leq
c_3x\log(x)+b_1$;

(c) For every $j$ we have $ord_{p_j}(F_{\ell_i})\geq c_4x(m-1)-b_2$.
\endlemma

The Lemma implies the theorem: indeed, by theorem
\ref{estimateinp2}, the $P_{\ell_i}$ verify the hypotheses of
proposition \ref{linearindependencecriterion2}.

\smallskip
\Proof (of lemma \ref{computationatq}) The only thing we have to
prove is (a); indeed (b) and (c) are consequence of proposition
\ref{globalsections}.

Let $\alpha$ be the order of $(P_1\wedge\dots\wedge P_m)$ at $q$.
We know that $\alpha\leq cx\epsilon+c_4$. By induction we have
that, if $h_1<\dots <h_m$,
$$\nabla_\partial(P_{h_1}\wedge\dots\wedge P_{h_m})=\sum_{s_1<\dots
<s_m}a_s\cdot(P_{s_1}\wedge\dots\wedge P_{s_m})$$ with $s_i\leq
h_m+1$ and suitable universal constants $a_s$. Consequently,
denoting by $\nabla^{\circ \alpha}_\partial(\cdot)$ the operator
$\nabla_\partial\circ\dots\nabla_\partial(\cdot)$ ($\alpha$
times), we have that
$$0\neq\nabla^{\circ \alpha}_\partial(P_1\wedge\dots\wedge
P_m)\vert_q=\sum_{s_1<\dots <s_m}a_s\cdot(P_{s_1}\wedge\dots\wedge
P_{s_m})\vert_q$$ with $s_m\leq m+cx+a$. Thus there exists
$\ell_1<l_2<\dots <\ell_m\leq m +cx+a$ such that
$$(P_{\ell_1}\wedge\dots\wedge P_{\ell_m})\vert_q\neq 0.$$
The conclusion follows.

\endsection

\section References\par

\

\item{[A1]} Andr\'e, Yves {\it S\'eries Gevrey de type arithm\'etique. I. Th\'eor${\rm\grave{e}}$mes de puret\'e et de dualit\'e}.
Ann. of Math. (2) 151 (2000), no. 2, 705--740.

\item{[A2]} Andr\'e, Yves {\it S\'eries Gevrey de type arithm\'etique. II. Transcendance sans
transcendance.}  Ann. of Math. (2) 151 (2000), no. 2, 741--756.

\item{[Ax]} Ax, James,{\it On Schanuel's conjectures}. Ann. of Math. (2) 93 1971 252--268.

\item{[Be1]} Bertrand, Daniel {\it Un th\'{e}or\`{e}me de Schneider-Lang sur certains domaines non simplement
connexes.} S\'eminaire Delange-Pisot-Poitou (16e ann\'ee:
1974/75), Th\'eorie des nombres, Fasc. 2, Exp. No. G18, 13 pp.
Secr\'etariat Math\'ematique, Paris, 1975.

\item{[Be2]} Bertrand, Daniel {\it Le th\'{e}or\`{e}me de
Siegel Shidlowsky revisit\'e}. Preprint.

\item{[Be3]}  Bertrand, Daniel {\it On Andr\'e's proof of the Siegel-Shidlovsky theorem.}
Colloque Franco -- Japonais: Th\'eorie des Nombres Transcendants
(Tokyo, 1998), 51--63

\item{[Beu]} Beukers, F. {\it A refined version of the Siegel-Shidlovskii theorem.} Ann. of Math. (2) 163 (2006), no. 1, 369--379.

\item{[Ga1]} Gasbarri, C. {\it Analytic subvarieties with many
rational points}, preprint available at

http://xxx.lanl.gov/abs/0811.3195

\item{[La]} Lang, Serge {\it Introduction to transcendental numbers.}
Addison-Wesley Publishing Co., Reading, Mass.-London-Don Mills,
Ont. 1966 vi+105 pp.

\item {[Ma]} Malgrange, B. {Sur les d\'eformations isomonodromiques. I. Singularit\'es r\'egulieres. II. Singularit\'es irr\'eguli\'eres.}
Mathematics and physics (Paris, 1979/1982), 427--438, Progr.
Math., 37, Birkh\"auser Boston, Boston, MA, 1983.

\item{[Vo]}  Vojta, Paul, {\it Diophantine approximations and value distribution theory.}
Lecture Notes in Mathematics, 1239. Springer-Verlag, Berlin, 1987.
x+132 pp.

\endsection

\bigskip

C.~Gasbarri:  Universit\'e de Strasbourg, IRMA,  7 rue Ren\'e Descartes
 
67084 Strasbourg -- France.

\

Email: gasbarri@math.u-strasbg.fr

\end